\documentclass{amsproc}

\usepackage{a4wide}
\usepackage{amsmath}
\usepackage{enumerate}

\usepackage{amsfonts, amssymb,amsthm,amsgen}

\newtheorem{theorem}{Theorem}[section]
\newtheorem{prop}{Proposition}[section]

\newtheorem{lem}{Lemma}[section]

\newtheorem{corol}{Corollary}[section]

\def\zr{\ltimes}
\def\ZL{\rightthreetimes}

\def\spa{\mathop\text{{\rm span}}\nolimits}

\def\pr{\mathop\text{\rm pr}\nolimits}

\def\Real{\mathbb{R}}
\def\Co{\mathbb{C}}

\def\g{\mathfrak{g}}
\def\h{\mathfrak{h}}
\def\so{\mathfrak{so}}

\def\sl{\mathfrak{sl}}
\def\gl{\mathfrak{gl}}
\def\su{\mathfrak{su}}

\def\simil{\mathfrak{sim}}

\def\u{\mathfrak{u}}
\def\m{\mathfrak{m}}

\def\f{\mathfrak{f}}
\def\k{\mathfrak{k}}

\def\m{\mathfrak{m}}
\def\z{\mathfrak{z}}
\def\R{\mathcal{R}}

\def\Simil{\mathop\text{\rm Sim}\nolimits}

\def\Un{\mathop\text{\rm U}\nolimits}
\def\SU{\mathop\text{\rm SU}\nolimits}

\def\Og{\mathop\text{\rm O}\nolimits}

\def\id{\mathop\text{\rm id}\nolimits}

\def\im{\mathop\text{\rm Im}\nolimits}
\def\re{\mathop\text{\rm Re}\nolimits}

\title{Holonomy classification of Lorentz-K\"ahler manifolds}

\author{Anton S. Galaev}\thanks{University of Hradec Kr\'alov\'e, Faculty of Science, Rokitansk\'eho 62, 500~03 Hradec Kr\'alov\'e,  Czech Republic\\
E-mail: anton.galaev(at)uhk.cz}

\begin{document}

\maketitle\vskip-50ex {\renewcommand{\abstractname}{Abstract}\begin{abstract}

The classification problem for holonomy of pseudo-Riemannian
manifolds is actual and open. In the present paper, holonomy
algebras of Lorentz-K\"ahler manifolds are classified. A simple
construction of a metric for each holonomy algebra is given. Complex
Walker coordinates are introduced and described using the potential.
Complex pp-waves are characterized in terms of the curvature,
holonomy and the potential. Classification of Lorentz-K\"ahler
symmetric spaces is reviewed.

\end{abstract}


{\bf 2010 Mathematical Subject Codes:} 53C29, 53C25, 53C35, 53C50,
53C55

\tableofcontents

\section{Introduction}

The holonomy group of a pseudo-Riemannian manifold $(M,g)$ is an
important invariant that gives rich information about the geometry
of $(M,g)$. This motivates the classification problem for holonomy
groups of pseudo-Riemannian manifolds. The problem is solved only
for connected holonomy groups of Riemannian and Lorentzian
manifolds. The case of Riemannian manifolds is classical and it has
many applications in differential geometry and theoretical physics,
see e.g. \cite{Besse,Bryant06,Joyce07}. The case of Lorentzian
manifolds took the attention of geometers and theoretical physicists
during the last two decades, see the reviews \cite{Baum12,ESI,UMS}
and the references therein. In the other signatures, only partial
results are known \cite{B13,B11,BB-I97,F-K,IRMA,DisB,G-T01,Ik2n}.
The main difference  between holonomy groups of Riemannian and
proper pseudo-Riemannian manifolds is that for Riemannian manifolds
all considerations may be reduced to the case of irreducible
holonomy groups, while in the case of proper pseudo-Riemannian
manifolds one should consider also holonomy groups preserving
degenerate subspaces. Berger \cite{Ber} classified connected
irreducible holonomy groups of pseudo-Riemannian manifolds of
arbitrary signature. Note that there is a more general
classification of connected irreducible holonomy groups of
torsion-free affine connections~\cite{M-Sch99}.

The general case of the signature $(2,N)$, i.e., the next signature
after the Lorentzian one, is quite complicated, this show the work
\cite{Ik2n}. In the present paper we give a complete classification
of the connected holonomy groups (equivalently, of holonomy
algebras) for pseudo-K\"ahler manifolds of complex index one; such
manifolds are called Lorentz-K\"ahler and they have holonomy
algebras contained in $\u(1,n+1)\subset\so(2,2n+2)$ (here $n+2$ is
the complex dimension of the manifold).

We do the same three steps that were previously done for the
holonomy of Lorentzian manifolds \cite{BB-I93,Leistner,Gmetrics}.
The Wu theorem allows to assume that the holonomy algebra is weakly
irreducible (this means that  it does not preserve any proper
non-degenerate subspace of the tangent space). Since the case of
irreducible holonomy algebras is solved by Berger, we assume that
the holonomy algebra is not irreducible, and then (with one
exception in dimension four) it preserves a complex isotropic line.
In Section \ref{secW-Ir} we classify  weakly irreducible subalgebras
in $\u(1,n+1)$. For that we consider the following geometrical
approach: we consider the induced action of the corresponding
connected Lie subgroup $G\subset\Un(1,n+1)$ on the boundary of the
complex hyperbolic space, which may be identified with the sphere.
This action induces an action of $G$ on the Heisenberg space and on
the Hermitian space $\Co^n$. Then $G$ preserves a real affine
subspace $L\subset\Co^n$, where it acts transitively as the group of
similarity transformations. This allows to find all groups $G$.
Using this, in Section \ref{secBerger} we  classify weakly
irreducible Berger subalgebras in $\u(1,n+1)$, i.e. weakly
irreducible subalgebras in $\u(1,n+1)$ having the same algebraic
properties as the holonomy algebras. Berger algebras play the role
of the candidates to the holonomy algebras.

Then appears the most important problem to check, which of these
candidates are really holonomy algebras. Note that in the case of
Riemannian manifolds the same problem was open for more than thirty
years, and finally it was solved by Bryant \cite{Bryant87}. In
Section \ref{secLKmetrics} we show that some of the obtained by us
Berger algebras cannot appear as the holonomy algebras (note that
all previously known Berger algebras were always holonomy algebras),
and we construct metrics realizing all the other Berger algebras as
holonomy algebras. We obtain the metrics by writing down explicitly
their potentials, which are relatively simple.

In Section \ref{secLKmetrics} we also introduce the complex
version of the Walker coordinates, which are frequently used in
Lorentzian geometry and Relativity \cite{WalkerBook}.

Then we consider several examples. In Section \ref{secLines} we
find the holonomy algebra of the space of oriented lines in
$\Real^3$, which admits the structure of the Lorentz-K\"ahler
manifold of real dimension four; this space is used in geometric
optics, see the survey \cite{G-K08} and the references therein.
There are generalizations of this space giving farther examples of
Lorentz-K\"ahler metrics \cite{G-K08,G-K08A,G-G10}.

The pp-wave metrics are very useful in theoretical physics. In
Section \ref{secpp-w} we consider their generalization to the
complex case. We give equivalent conditions in terms of the
curvature, holonomy and the potential for a complex Walker metric to
be a complex pp-wave. Complex pp-waves were studied in \cite{L-L13},
where appear also examples of Lorentz-K\"ahler metrics and their
holonomy. Complex pp-waves are used also in the physical literature,
e.g. \cite{Chudecki,Plebanski}.

Simply connected pseudo-Riemannian symmetric spaces of index 2 were
classified in \cite{K-O09}, and the case of Lorentz-K\"ahler
manifolds was considered separately. In Section \ref{secSym} we
reformulate the last result in terms of the curvature and holonomy.

Let us mention the works and situations, where Lorentz-K\"ahler
manifolds and their holonomy algebras appear. There are various
examples of invariant Lorentz-K\"ahler metrics on Lie groups and
homogeneous spaces, e.g. \cite{B-Z11,C-F13,C-O12,O06,Yamada}. In
some cases the holonomy groups are computed. Lorentz-K\"ahler
metrics and their holonomy are related to  conformal geometry,
conformal holonomy, and CR-structures through the ambient metric
construction, see e.g. \cite{B-J10,C-G10,F-H03,Leitner08}, in
particular, to the conformal analog of Calabi-Yau manifolds
\cite{Baum10}.

Completing the introduction, we note that  the problem solved in the
present paper was  considered in the unpublished PhD thesis
\cite{DisB} of the author. The classification result from
\cite{DisB}  with a sketch of the proof was given in the survey
\cite{IRMA}, see remark in Section \ref{secMainTh}. The exposition
of the present paper is much simpler than that in \cite{DisB}, in
particular, the most important is the new construction of metrics
given in the present paper.

\vskip1cm

{\bf Acknowledgements.} The author is grateful to Helga Baum and
Dmitri V. Alekseevsky for useful discussions and suggestions. This
work was supported by the Specific Research Project of Faculty of
Science, University of Hradec Kr\'alov\'e [No. 2101, 2017] and by
the grant  no. 18-00496S of the Czech Science Foundation.

\section{Preliminaries}

The theory of the holonomy groups of pseudo-Riemannian manifolds
can be found e.g. in \cite{Besse,Joyce07,IRMA}. Let $(M,g)$ be a
pseudo-Riemannian manifold of signature $(r,s)$. {\it The holonomy
group} $G_x$ of $(M,g)$ at a point $x\in M$ is the  Lie group that
consists of the pseudo-orthogonal  transformations of the tangent
space given by the parallel transports along all piecewise smooth
loops at the point $x$. This group may be identified with a Lie
subgroup of the pseudo-orthogonal group $\Og(r,s)$. The
corresponding subalgebra of the pseudo-orthogonal Lie algebra
$\so(r,s)$ is called {\it the holonomy algebra} and it determines
the holonomy group if the manifold is simply connected. {\it The
Ambrose-Singer Theorem} states that the holonomy algebra is
spanned by the endomorphisms
 $\tau^{-1}_\gamma\circ R_y(\tau_\gamma X,\tau_\gamma
Y)\circ \tau_\gamma$ of $T_xM$, where  $\gamma$ is a piecewise
smooth curve starting at the point $x$ with an end-point $y \in
M$, and $X,Y\in T_xM$. This implies that the holonomy algebra
satisfies a strong algebraic condition being spanned by images of
the algebraic curvature tensors;  algebras satisfying this
property are called Berger algebras (see below), and they are
candidates to the holonomy algebras.

 The {\it fundamental principle for holonomy groups}
states that there exists a one-to-one correspondents between
parallel tensor fields $T$  on $(M,g)$ and tensors $T_0$ of the
same type at $x$ preserved by the tensor extension of the
representation of the holonomy group.

A subgroup (resp. subalgebra) of the pseudo-orthogonal Lie group
(resp. algebra) is called  {\it weakly irreducible} if it does not
preserve any non-degenerate proper subspace of the
pseudo-Euclidean space). The Wu decomposition theorem implies that
if the holonomy algebra of a pseudo-Riemannian manifold is not
weakly irreducible, then the manifold is locally decomposable, and
the holonomy algebra is the direct sum of weakly irreducible
holonomy algebras. This allows us to assume that the holonomy
algebra is weakly irreducible. Irreducible holonomy algebras of
pseudo-Riemannian manifolds were classified by Berger \cite{Ber},
so we are left with weakly irreducible not irreducible holonomy
algebras.

{\it A pseudo-K\"ahler manifold} is a pseudo-Riemannian manifold
$(M,g)$ with a $g$-orthogonal complex structure $J$ that is
parallel. By the fundamental principle, the equivalent condition
is that the holonomy group of $(M,g)$ is contained in the
pseudo-unitary group $\Un(\frac{r}{2},\frac{s}{2})$ (the numbers
$r$ and $s$ must be even). Together with $J$, $M$ becomes a
complex manifold; the tensor fields $g$ and $J$ define a
pseudo-K\"ahler metric $h$ on $M$ of complex signature
$(\frac{r}{2},\frac{s}{2})$.
 We are interested in holonomy groups of pseudo-K\"ahler
manifold of real signature $(2,2N)$, i.e. of complex signature
$(1,N)$. Such manifolds are called {\it Lorentz-K\"ahler manifolds}.
Let $N=n+1$, $n\geq 0$. The tangent space of a Lorentz-K\"ahler
manifold may be identified with the pseudo-Hermitian space
$\Co^{1,n+1}$ with a pseudo-Hermitian metric $h$. From the Berger
classification it follows that the only irreducible holonomy
algebras of Lorentz-K\"ahler  manifolds are $\u(1,n+1)$ and
$\su(1,n+1)$. In fact, these Lie algebras exhaust irreducible
subalgebras in $\u(1,n+1)$ \cite{D-L11}.

Fix a Witt basis $p,e_1,\dots,e_n,q$ of $\Co^{1,n+1}$, with respect
to such a basis the non-zero values of the metric $h$ are the
following: $h(p,q)=h(e_j,e_j)=1$. In particular, the vectors $p$ and
$q$ are isotropic, and the vectors $e_1,\dots,e_n$ form a basis of
the Hermitian space, which we denote by $\Co^n$. We will be
interested in the subalgebra $\u(1,n+1)_{\Co p}$ of $\u(1,n+1)$
preserving the complex isotropic line $\Co p$. In the matrix form we
have
$$\u(1,n+1)_{\Co p}=\left.\left\{\left(\begin{array}{ccc}
a&-\bar Z^t&ic\\
0&A&Z\\
0&0&-\bar a\end{array}\right)\right|\begin{array}{c}a\in\Co,\,A\in\u(n),\\
Z\in\Co^{n}\\ c\in\Real\end{array}\right\}.$$ We will denote the
above matrix by the 4-tuple $(a,A,Z,c)$. For the Lie brackets we
have
$$[(a,A,0,0),(b,B,Z,c)]=(0,[A,B]_{\u(n)},\bar a Z+AZ,2ic\re a),\quad
[(0,0,Z,0),(0,0,V,c)]=(0,0,0,2i\im h(Z,V)).$$ We obtain the
decomposition $$\u(1,n+1)_{\Co p}=(\Co\oplus\u(n))\zr(\Co^n\zr
i\Real).$$ Note that the subalgebra $\Co^n\zr i\Real\subset
\u(1,n+1)_{\Co p}$ is isomorphic to the Heisenberg Lie algebra,
and the subalgebra $\su(1,n+1)_{\Co p}\subset \u(1,n+1)_{\Co p}$
is isomorphic to the Lie algebra of the similarity transformations
of the Heisenberg space. We will denote the complex structure on
$\Co^{1,n+1}$ by $J$.

We will consider real vector subspaces of $\Co^n$ of the form
$$L=\Co^m\oplus L_0,$$
where $0\leq m \leq n$, we fix an $h$-orthogonal  decomposition
$\Co^n=\Co^{m}\oplus\Co^{n-m}$, and $L_0\subset\Co^{n-m}$ is a
real form, i.e. $iL_0\cap L_0=0$, $L_0\oplus iL_0=\Co^{n-m}$. An
example of $L_0$ is
$$\Real^{n-m}=\spa_\Real\{e_{m+1},\dots,e_n\}.$$ Note that in
general $L_0$ is different from this one, since $L_0$ may not
contain an $h$-orthonormal basis of $\Co^{n-m}$.

For the case $n=0$ let us consider the  subalgebra
$\g_0\subset\su(1,1)$ defined in the following way. Consider the
Witt basis $p_1,p_2,q_1,q_2$ in $\Real^{2,2}$ and define the
pseudo-Hermitian structure $J$ such that $Jp_1=q_2$, $Jp_2=-q_1$.
Let
$$\g_0=\left.\left\{\left(\begin{array}{cc}
A&0\\
0&-A^t\end{array}\right)\right|A\in\sl(2,\Real)\right\}\subset\so(2,2).$$
The Lie algebra $\g_0$ commutes with $J$, hence it is contained in
$\u(1,1)$. Moreover, $\g_0\subset\su(1,1)$. The Lie algebra $\g_0$
preserves the vector subspace $\spa_\Real\{p_1,p_2\}$, which is
not $J$-invariant.

\section{The Classification Theorem}\label{secMainTh}

\begin{theorem}\label{ThMain}

{\bf 1)} A subalgebra $\g\subset\u(1,1)$ is the weakly irreducible
not irreducible holonomy algebra of a Lorentz-K\"ahler manifold of
complex dimension 2 if and only if $\g$ is conjugated to $\g_0$,
or $\g$ is conjugated to one of the following subalgebras of
$\u(1,1)_{\Co p}$:
\begin{itemize}
\item $\g_1=\u(1,1)_{\Co p}$;

\item $\g_2=\left.\left\{\left(\begin{array}{cc}
a&0\\
0&-\bar a\end{array}\right)\right|a\in\Co\right\};$

\item
$\g_3^\gamma=\left.\left\{\left(\begin{array}{cc}
r\gamma&ic\\
0&-r\bar \gamma\end{array}\right)\right|r,c\in\Real\right\},$
 where
$\gamma\in\Co$ is a fixed number.
\end{itemize}

{\bf 2)} Let $n\geq 1$. Then a subalgebra $\g\subset\u(1,n+1)$ is
the weakly irreducible not irreducible holonomy algebra of a
Lorentz-K\"ahler manifold of complex dimension $n+2$ if and only if
$\g$ is conjugated to one of the following subalgebras of
$\u(1,n+1)_{\Co p}$:
\begin{itemize}

\item $\g^{\k}=\k\zr(\Co^n\zr i\Real)$\\
$=\left.\left\{\left(\begin{array}{ccc}
a&-\bar Z^t&ic\\
0&A&Z\\
0&0&-\bar a\end{array}\right)\right|\begin{array}{c}a+A\in\k,\\
Z\in\Co^{n}\\ c\in\Real\end{array}\right\},$\\ where
$$\k\subset\Co\oplus\u(n)$$ is an
arbitrary subalgebra;

\item $\g^{\k,J,L}=\k\zr(L\zr i\Real)$\\
$=\left.\left\{\left(\begin{array}{cccc}
a_2i&-\bar Z^t&-\bar X^t &ic\\
0&A&0&Z\\
0&0& a_2iE_{n-m}&X\\
0&0&0&a_2i\end{array}\right)\right|\begin{array}{c}a_2(i+i\id_{\Co^{n-m}})+A\in\k,\\
Z\in\Co^{m},\,X\in \Real^{n-m},\\ c\in\Real\end{array}\right\},$\\
where $L=\Co^{m}\oplus \Real^{n-m}$, $0\leq m<n$,
$$\k\subset \Real J\oplus\u(m)$$ is an
arbitrary subalgebra not contained in $\u(m)$;

\item $\g^{\k,L}=\k\zr(L\zr i\Real)$\\
$=\left.\left\{\left(\begin{array}{cccc}
0&-\bar Z^t&-\bar X^t &ic\\
0&A&0&Z\\
0&0& 0&X\\
0&0&0&0\end{array}\right)\right|\begin{array}{c}A\in\k,\\
Z\in\Co^{m},\,X\in L_0,\\ c\in\Real\end{array}\right\},$\\ where
$L=\Co^{m}\oplus L_0$, $0\leq m<n$, $L_0\subset\Co^{n-m}$ is a
real form,
$$\k\subset \u(m)$$ is an
arbitrary subalgebra;

\item $\g^{\k_0,\psi}=\k_0\oplus\{(0,\psi(X),X,0)|X\in \Co^{m-r}\oplus L_0\}\zr(\Co^r\zr i\Real)$
$$=\left.\left\{\left(\begin{array}{cccc}
0&-\bar Z^t&-\bar X^t &ic\\
0&A+\psi(X)&0&Z\\
0&0& 0&X\\
0&0&0&0\end{array}\right)\right|\begin{array}{c}A\in\k_0,\\
Z\in\Co^{r},\,X\in \Co^{m-r}\oplus L_0,\\
c\in\Real\end{array}\right\},$$ where $1\leq r\leq m\leq n$,
$L=\Co^{m}\oplus L_0$,  $L_0\subset\Co^{n-m}$ is a real form,
 we
consider the decomposition
$$\Co^m=\Co^{r}\oplus\Co^{m-r},$$ $\k_0\subset\u(r)$ is an arbitrary
subalgebra,
$$\psi:\Co^{m-r}\oplus L_0\to \u(r)$$ is a non-zero linear map such
that the image $\psi(\Co^{m-r}\oplus L_0)\subset\u(r)$ is
commutative, it commutes with $\k_0$ and has trivial intersection
with $\k_0$.
\end{itemize}

\end{theorem}

Recall that, by definition, a special pseudo-K\"ahler manifold (or
Calabi-Yau manifold of indefinite signature) is just a Ricci-flat
pseudo-K\"ahler manifold. Special pseudo-K\"ahler manifolds are
precisely pseudo-K\"ahler manifolds with holonomy algebras contained
in the special unitary Lie algebra. We imminently get the following
corollary.

\begin{corol} Weakly irreducible not irreducible holonomy algebras
of Ricci-flat Lorentz-K\"ahler manifolds are exhausted by the
following algebras:

 $\g_0$;

 $\g_3^\gamma$ with $\gamma\in\Real$;

$\g^\k$ with
$\k\subset\Real\oplus\Real(ni-2\id_{\Co^n})\oplus\su(n)$,

$\g^{\k,J,L}$ with $\k\subset\Real\left(mi-(2+n-m)i\id_{\Co^m}+
mi\id_{\Co^{n-m}}\right)\oplus\su(m)$;

 $\g^{\k,L}$ with
$\k\subset\su(m)$;

 $\g^{\k_0,\psi}$ with $\k_0\oplus\psi(\Co^{m-r}\oplus
L_0)\subset\su(r)$.
\end{corol}

{\bf Proof of Theorem \ref{ThMain}.}

Suppose that $\g\subset\u(1,n+1)$ is the weakly irreducible not
irreducible holonomy algebra of a Lorentz-K\"ahler manifold. Then
$\g$ preserves a proper real subspace
$W\subset\Real^{2,2n+2}=\Co^{1,n+1}$. Consequently $\g$ preserves
the orthogonal complement $W^\bot\subset\Real^{2,2n+2}$ and the
intersection $W\cap W^\bot$, which is isotropic and must be of
real dimension 1 or 2. Let $W_1=(W\cap W^\bot)^\bot$.

{\bf Suppose that $n=0$}. If $W_1\cap JW_1\neq 0$, then $\g$ is
conjugated to a subalgebra of $\u(1,1)_{\Co p}$. It is not hard to
write down all subalgebras of $\u(1,1)_{\Co p}$ and decide which of
them are  weakly irreducible Berger subalgebras in $\u(1,1)_{\Co
p}$, i.e. weakly irreducible subalgebras in $\u(1,1)_{\Co p}$ having
the same algebraic properties as the holonomy algebras. In Section
\ref{secLKmetrics} we construct metrics with each of the obtained
algebras being the holonomy algebra. If $W_1\cap JW_1=0,$ then it is
not hard to see that $\g$ is contained in $\g_0$. Holonomy algebras
of pseudo-Riemannian manifolds of signature $(2,2)$ are well studied
\cite{BB-I97,G-T01}, by this reason we do not pay much attention to
the case $n=0$.

{\bf Suppose that $n>0$}. Then $W_2=W_1\cap JW_1$ is non-trivial,
$\g$- and $J$-invariant, and degenerate. Thus, $\g$ preserves the
complex isotropic line $W_2\cap W_2^\bot$. Hence $\g$ is
conjugated to a subalgebra of $\u(1,n+1)_{\Co p}.$

The next three sections will be dedicated to the rest of the proof
of Theorem \ref{ThMain}. In Section \ref{secW-Ir} we will describe
weakly irreducible subalgebras in $\u(1,n+1)_{\Co p}$. Using this,
in Section \ref{secBerger} we  classify weakly irreducible Berger
subalgebras in $\u(1,n+1)_{\Co p}$. In Section \ref{secLKmetrics} we
show that some of the obtained Berger algebras cannot appear as the
holonomy algebras, and we construct metrics realizing all the other
Berger algebras as holonomy algebras. This will prove Theorem
\ref{ThMain}. \qed

\medskip

{\bf Remark.} The first attempt to classify holonomy algebras of
Lorentz-K\"ahler manifolds was done in the unpublished PhD thesis
\cite{DisB} of the author. The exposition of \cite{DisB} is very
complicated compared to the exposition of the present paper, e.g.,
the metrics constructed in \cite{DisB} were given in real
coordinates and in a very involved way. In \cite{DisB}, it was
wrongly assumed that any real form $L\subset\Co^n$ contains an
orthonormal basis for the Hermitian form on $\Co^n$, i.e., the
endomorphism $\theta$ defined below was assumed to be zero. The
statement of Theorem \ref{ThnonEx} given below was not discovered in
\cite{DisB}. The classification result from \cite{DisB} with a
sketch of the proof was given in the survey \cite{IRMA} (Theorem
4.2). Let us compare \cite[Th. 4.2]{IRMA} with the above Theorem
\ref{ThMain} in the case $n\geq 1$.  The results of the present
paper imply that each of the first five algebras $\h$ from \cite[Th.
4.2, part 2]{IRMA} must satisfy the condition: if
$\h\not\subset(\Real J\oplus\u(n))\zr(\Co^n\zr i\Real)$, then the
number $m$ associated to $\h$ is equal to $n$ (this is the
consequence of Theorem \ref{ThnonEx} given below). Under that
condition these algebras exhaust the first two algebras from part 2
of Theorem \ref{ThMain}. The third algebra from part 2 of Theorem
\ref{ThMain} that satisfies the condition $L_0\neq \Real^{n-m}$
(i.e., with $\theta\neq 0$) is missing in \cite[Th. 4.2]{IRMA}. If
$L_0= \Real^{n-m}$ (i.e., $\theta= 0$) then this is the third
algebra from \cite[Th. 4.2, part 2]{IRMA} with $\varphi=\phi=0$.
Likewise, the last two algebras from \cite[Th. 4.2]{IRMA} correspond
to the last algebra from Theorem \ref{ThMain} with
$L_0=\Real^{n-m}$.

\section{Weakly irreducible subalgebras in $\u(1,n+1)$}\label{secW-Ir}

Let $\Un(1,n+1)_{\Co p}$ be the Lie subgroup of $\Un(1,n+1)$
preserving the line $\Co p$. This group is generated by the
elements of the form:

$$\left (\begin{array}{ccc} e^a & 0 & 0\\ 0 & A & 0\\ 0 & 0 &
e^{-\bar a}\end{array}\right),\quad \left( \begin{array}{ccc} 1 &
- \bar Z^t & -
\frac{1}{2} \bar Z^tZ+ic\\ 0 & E_n & Z\\ 0 & 0& 1 \\
\end{array} \right),$$
where $a\in\Co$, $A\in \Un(n)$, $Z\in\Co^n$, and $c\in\Real$.

Consider the  group $$\Simil \Co^n=(\Real^*\cdot U(n))\ZL \Co^n$$
 of similarity transformations of the
Hermitian space $\Co^n$ and the group homomorphism
$$\Gamma:\Un(1,n+1)_{\Co p}\to \Simil \Co^n $$
that sends the above elements of $\Un(1,n+1)_{\Co p}$ to $e^{\bar
a}\cdot A\in \Real^*\cdot U(n)$, and $Z\in\Co^n$, respectively.

The homomorphism $\Gamma$ can be constructed  in the following
geometrical way. The group $\Un(1,n+1)$ acts in the natural way on
the boundary $\partial \mathbf{H}^{n+1}$ of the complex hyperbolic
space, which consists of complex isotropic lines in $\Co^{1,n+1}$
\cite{Goldman}. The boundary may be identified with the sphere
$S^{2n+1}$. Since the group $\Un(1,n+1)_{\Co p}$ preserves the
point $\{\Co p\}\in
\partial \mathbf{H}^{n+1}$, it acts on the sphere with one removed
point; this space may be identified with the Heisenberg space
$$\Co^n\oplus\Real,$$
and the action there  is given by the Heisenberg similarity
transformations. Then we may consider the induced action on
$\Co^n$, which will give us the homomorphism $\Gamma$.

The differential
$$\Gamma':\u(1,n+1)_{\Co
p}\to\simil \Co^n=(\Real\oplus\u(n))\zr\Co^n$$ is given by
$$(a,A,Z,c)\mapsto (\re a,-i\im a\id_{\Co^n}+A,Z).$$
The kernel of $\Gamma'$ equals to $\Real J\oplus i\Real\subset
\u(1,n+1)_{\Co p}$.

Suppose that $G\subset\Un(1,n+1)_{\Co p}$ is a weakly irreducible
Lie subgroup in the sense that it does not preserve any
non-degenerate complex subspace in $\Co^{1,n+1}$.

\begin{prop}\label{propnon-compl}
The subgroup $\Gamma(G)\subset \Simil \Co^n$ does not preserve any
proper complex affine subspace of $\Co^n$. Consequently,
 if $\Gamma(G)\subset \Simil \Co^n$ preserves a proper
real affine subspace $L\subset \Co^n$, then the minimal complex
affine subspace of $\Co^n$ containing $L$ is $\Co^n$.
\end{prop}

{\it Proof.}  First we prove that $\Gamma(G)\subset \Simil \Co^n$
does not preserve any proper complex  vector subspace of $\Co^n$.
Suppose that $\Gamma(G)$ preserves a proper complex vector
subspace $L\subset \Co^n$. Then, $$\Gamma(G)\subset (\Real^*\cdot
U(L)\cdot U(L^{\bot}))\ZL L,$$ where $L^{\bot}$ is the orthogonal
complement to $L$ in $\Co^n$. We see that the group $G$ preserves
the proper non-degenerate vector subspace $L^{\bot}\subset
\Co^{1,n+1}$.

Suppose that $\Gamma(G)\subset \Simil \Co^n$ preserves a proper
complex affine subspace $L\subset \Co^n$. The subgroup  $G\subset
\Un(1,n+1)_{\Co p}$ is conjugated to a subgroup
$G_1\subset\Un(1,n+1)_{\Co p}$ such that $\Gamma(G_1)\subset \Simil
\Co^n$ preserves a proper complex vector subspace corresponding to
the affine subspace $L\subset \Co^n$. This gives a contradiction,
since it is clear that $G_1\subset\Un(1,n+1)_{\Co p}$ is weakly
irreducible.
 This proves the proposition. \qed

Let $\g\subset\u(1,n+1)_{\Co p}$ be a weakly irreducible
subalgebra. Let $G\subset\Un(1,n+1)_{\Co p}$ be the corresponding
connected Lie subgroup. Suppose that  $\Gamma(G)$ preserves a real
affine subspace $L\subset\Co^n$. Considering, if necessary,
another group $G_1$ conjugated to $G$, we may assume that
$L\subset\Co^n$ is a real vector subspace. We may assume also that
$L$ is minimal in the sense that $\Gamma(G)$ does not preserve any
proper affine subspace in $L$. The induced action of $\Gamma(G)$
on the Euclidean space $L$ is by similarity transformations. Since
$\Gamma(G)$ does not preserve any proper affine subspace of $L$,
it acts transitively on $L$ \cite{A-V-S}. Thus,
$$\Gamma(G)\subset (\Og(L^\bot)\cdot \Simil L)\cap\Simil\Co^n=
(\Real^*\cdot ((\Og(L^\bot)\cdot\Og(L))\cap\Un(n)))\ZL L,$$ where
the orthogonal complement is taken with respect to the Euclidean
metric $\re h$, and
$$\Simil L=(\Real^*\cdot \Og(L))\ZL L$$
is the  Lie group of similarity transformations of the Euclidean
space $L$. Moreover, the projection of $\Gamma(G)$ on $\Simil L$ is
a transitive group of similarity transformations. The Lie algebras
of the transitive Lie groups of similarity transformations of an
Euclidean space $E$ are exhausted by the following \cite{A-V-S}:
\begin{itemize}
\item $\f\zr L$,
\item $(\f_0\oplus\{\psi(X)+X|X\in U\})\zr W$,
\end{itemize}
where $\f\subset\Real\oplus\so(L)$ is a subalgebra;  for the
second algebra, there exists an orthogonal decomposition
$L=W\oplus U$, $\f_0\subset\so(W)$ is a subalgebra,
$\psi:U\to\so(W)$ is a non-zero linear map such that
$\f_0\cap\psi(U)=0$ and $\psi(U)$ is commutative and it commutes
with $\f_0$.

Let us describe the subspace $L\subset\Co^n$. It is clear that
$iL\cap L\subset \Co^n$ is a Hermitian subspace, let us denote it
by $\Co^m$. Let $\Co^{n-m}$ denote the orthogonal complement to
$\Co^m$. We obtain the decomposition
$$\Co^n=\Co^m\oplus\Co^{n-m}.$$
Let $L_0\subset L$ be the orthogonal complement to $\Co^m\subset
L$ with respect to $\re h$. It is clear that $iL_0\cap L_0=0$,
i.e. $L_0\subset \Co^{n-m}$ is a real form. We obtain the
decomposition $$L=\Co^m\oplus L_0.$$ Let $e_1,...,e_m$ be an
$h$-orthonormal basis of $\Co^m$. Consider a basis
$f_{m+1},\dots,f_n$ of the real vector space $L_0$. We may assume
that this basis is orthonormal with respect to $\re h$. There
exists a skew-symmetric real matrix
$\omega=(\omega_{jk})_{j,k=m+1}^n$ such that
$$h(f_j,f_k)=\delta_{jk}+i\omega_{jk},\quad j,k=m+1,\dots,n.$$
It is known, that the basis $f_{m+1},\dots,f_n$ can be chosen in
such a way that
$$\omega=\textrm{diag}\left(\left(\begin{array}{cc}
0&-\lambda_1\\\lambda_1&0\end{array}\right),\dots,\left(\begin{array}{cc}
0&-\lambda_s\\\lambda_s&0\end{array}\right),0,\dots,0\right)$$ for
some real numbers  $\lambda_k$. Since $h$ is positive definite, it
holds $|\lambda_k|<1$.

Let $\tau:\Co^{n-m}\to\Co^{n-m}$ be the anti-linear involution
defining the real form $L_0\subset \Co^{n-m}$. It holds
$$\tau(f_k)=f_k,\quad \tau(if_k)=-if_k,\quad k=m+1,\dots,n.$$
Define the endomorphism $\theta$ of $L_0$ by the equation
$${\rm Re}h(\theta X, Y)={\rm Im} h(X,Y),\quad X,Y\in L_0$$
and extend it to a $\Co$-linear endomorphism of $\Co^{n-m}$.
Clearly, $\theta\in\u(n-m)$. Note that $\theta=0$ if and only if
$L_0$  contains an $h$-orthonormal basis of $\Co^{n-m}$.

For $k=1,\dots,s$, let
$$e_{m+2k-1}=\frac{\sqrt{2}}{2\sqrt{1-\lambda_k}}(f_{m+2k-1}+if_{m+2k}),\quad
e_{m+2k}=\frac{\sqrt{2}}{2\sqrt{1+\lambda_k}}(f_{m+2k}+if_{m+2k-1}).$$
Then it holds $$\theta(e_{m+2k-1})=-i\lambda_ke_{m+2k-1},\quad
\theta(e_{m+2k})=i\lambda_ke_{m+2k}$$ and
$$\tau(e_{m+2k-1})=-i\frac{\sqrt{1+\lambda_k}}{\sqrt{1-\lambda_k}}e_{m+2k},\quad
\tau(e_{m+2k})=-i\frac{\sqrt{1-\lambda_k}}{\sqrt{1+\lambda_k}}e_{m+2k-1}.$$

Note that
$$(\so(L^\bot)\oplus\so(L))\cap\u(n)=\u(m)\oplus\z_{\u(n-m)}\tau,$$
where $\z_{\u(n-m)}\tau$ denotes the ideal in $\u(n-m)$ consisting
of elements commuting with $\tau$.

Suppose that the projection of $\Gamma'(\g)$ to $\simil(L)$ is of
the form $(\f_0\oplus\{\psi(X)+X|X\in U\})\zr W$. Since
$\f=\f_0\oplus\psi(U)$ annihilates $U$, and it is contained in
$\u(n)$, it annihilates also $iU$. Consequently,  $\f$ is contained
in $\u((U+iU)^\bot)$. We denote this space by $\Co^r$. It is clear
that $\Co^r\subset L$. This implies that $\Co^r\subset\Co^m$. Let
$\Co^{m-r}$ be the orthogonal complement to $\Co^r$ in $\Co^m$. We
get the decomposition
$$L=\Co^r\oplus\Co^{m-r}\oplus L_0$$
and we extend the map $\psi$ to $\Co^{m-r}\oplus L_0$. Thus we may
assume that $W=\Co^r$, and $U=\Co^{m-r}\oplus L_0$.

We conclude that if  $\g\subset\u(1,n+1)_{\Co p}$ is a weakly
irreducible subalgebra with the associated subspace $L=\Co^m\oplus
L_0\subset\Co^n$, then $\Gamma'(\g)$ is one of the following:
\begin{itemize}
\item $\f\zr L$,
\item $(\f_0\oplus\{\psi(X)+X|X\in \Co^{m-r}\oplus L_0\})\zr \Co^r$,
\end{itemize}
where $\f\subset\Real\oplus\u(m)\oplus\z_{\u(n-m)}\tau$ is a
subalgebra; $\f_0\subset\u(r)$, $\psi:\Co^{m-r}\oplus L_0\to\u(r)$
is a non-zero linear map such that $\f_0\cap\psi(\Co^{m-r}\oplus
L_0)=0$ and $\psi(\Co^{m-r}\oplus L_0)$ is commutative and it
commutes with $\f_0$.

Let $\g_0$ be the projection of $\g$ to $\su(1,n+1)$. Note that if
$m>0$, then from the structure of the Lie brackets in
$\u(1,n+1)_{\Co p}$ it follows that $\g_0$ contains the ideal
$i\Real\subset\u(1,n+1)_{\Co p}$. Even more, below we will see
that all holonomy algebras $\g\subset\u(1,n+1)_{\Co p}$ contain
this ideal, thus it is enough to consider only such subalgebras.
Then, since the kernel of $\Gamma'|_{\su(1,n+1)_{\Co p}}$
coincides with $i\Real$, we immediately find that
$$\g_0=(\Gamma'|_{\su(1,n+1)_{\Co p}})^{-1}(\Gamma'(\g)).$$
Note that $$(\Gamma'|_{\su(1,n+1)_{\Co p}})^{-1}(i\id_{\Co^m})=
\frac{i}{n+2}\left(\begin{array}{cccc} -m & 0 &
0&0\\0&n+2-m&0&0\\0&0&-m&0\\0&0&0&-m\end{array}\right)+i\Real.$$
Now, for the Lie algebra $\g$ we have the following there
possibilities:
$$\g=\g_0,\quad\g=\g_0\oplus\Real J,\quad
\g=\{A+\zeta(A)J|A\in\g_0\},$$ where $\zeta:\g_0\to\Real$ is a
linear map zero on the commutant $\g_0'=[\g_0,\g_0]$.

Thus we conclude that if $\g\subset\u(1,n+1)_{\Co p}$ is a weakly
irreducible subalgebra containing the ideal $i\Real\subset
\u(1,n+1)_{\Co p}$, then $\g$ is one of the following:
\begin{itemize}
\item $\f\zr (L\zr i\Real)$,
\item $(\f_0\oplus\{\psi(X)+X|X\in \Co^{m-r}\oplus L_0\})\zr (\Co^r\zr i\Real)$,
\end{itemize}
where $\f\subset\Real\oplus\Real
J\oplus\u(m)\oplus\z_{\u(n-m)}\tau$ is a subalgebra;
$\f_0\subset\u(r)\oplus\Real J$, $\psi:\Co^{m-r}\oplus
L_0\to\u(r)\oplus\Real J$ is a non-zero linear map such that
$\f_0\cap\psi(\Co^{m-r}\oplus L_0)=0$ and $\psi(\Co^{m-r}\oplus
L_0)$ is commutative and it commutes with $\f_0$.

\section{Algebraic curvature tensors and Berger
algebras}\label{secBerger}

 Let $(M,g,J)$ be a Lorentz-K\"ahler manifold of complex dimension $n+2\geq 2$, i.e. $g$
 is a pseudo-Riemannian metric
 of signature $(2,2n+2)$ on $M$, and $J$ is a parallel $g$-orthogonal complex structure.
Consider the corresponding Lorentz-K\"ahler metric
$$h=g+ig\circ J.$$
 Fix a point $x\in
 M$.  Let $\g\subset\u(T_xM,h_x)$ be the holonomy algebra of the manifold $(M,g)$. We
 identify the  tangent space $T_xM$ with the Lorentz-Hermitian
 space $\Co^{1,n+1}$.

 Consider the complexified tangent bundle $T^\Co M$  and the
 standard decomposition
 $$T^\Co M=T^{1,0}M\oplus T^{0,1}M$$
 into the direct sum of eigenspaces of the $\Co$-linear extension
 of $J$ corresponding to the eigenvalues $\pm i$.
 Recall that sections of $T^{1,0}M$ and $T^{0,1}M$ are of the form
 $$X-iJX,\quad X+iJX,$$ where $X$ is a vector field on $M$. In
 particular, the bundle $T^{1,0}M$ is identified with $TM$.
The metric $g$ allows to identify the bundle $T^{0,1}M$ with the
dual bundle to $T^{1,0}M$. Moreover, let $X,Y$ denote both vector
fields in $TM$ and the corresponding vector fields in $T^{1,0}M$,
then it holds \begin{equation}\label{g=4h} g(X,\bar
Y)=4h(X,Y),\end{equation} where $g$ denotes the $\Co$-bilinear
extension to $T^\Co M$ of the initial metric $g$.

Fix a basis $p,e_1,\dots,e_n,q$ in $T_xM$ as above. Denote by the
same symbols $p,e_1,\dots,e_n,q$ the basis
$$\frac{1}{2}(p-iJp),\frac{1}{2}(e_1-iJe_1),\dots,\frac{1}{2}(e_n-iJe_n),\frac{1}{2}(q-iJq)\in T^{1,0}_xM.$$
Then the vectors $\bar p,\bar e_1,\dots,\bar e_n,\bar q$ form a
basis of $T_x^{0,1}M$. Suppose that $\g\subset\u(1,n+1)_{\Co p}$.
The holonomy algebra of the complex connection $\nabla$ coincides
with the complexification $\g^\Co=\g\otimes\Co$ of $\g$. The
complexification of $\u(1,n+1)_{\Co p}$ consists of the matrices
of the form
$$\left(\begin{array}{cc}T_1&0\\0&T_2\end{array}\right),$$
where
$$T_1=\left(\begin{array}{ccc}a&\bar W^t&c\\0&A&Z\\0&0&b\end{array}\right), \quad
T_2=-\left(\begin{array}{ccc}b&\bar
Z^t&c\\0&A^t&W\\0&0&a\end{array}\right)$$ with $a,b,c\in\Co$,
$Z,W\in \Co^n$ and  $A\in\gl(n,\Co)$. Denote by $\g^{1,0}$ the
projection of $\g^\Co$ to $\gl(T^{1,0}_xM)$. It is clear that this
Lie algebra determines $\g^\Co$. This Lie algebra allows us also
to determine $\g$, namely, consider the following anti-linear
involution of $\g^{1,0}$:
$$\sigma:\left(\begin{array}{ccc}a&\bar
W^t&c\\0&A&Z\\0&0&b\end{array}\right)\mapsto
-\left(\begin{array}{ccc}\bar b&\bar Z^t&\bar c\\0&\bar
A^t&W\\0&0&\bar a\end{array}\right).$$ Then $\g$ coincides with
the set of fixed points of $\sigma$.

The symbol $\nabla$ will  denote both the Levi-Civita connection on
$(M,g)$ and the corresponding complex connection on $T^\Co M$. We
use the same convection for $R$. It is well known that $R$ satisfies
the conditions
$$ R(JX,JY)=R(X,Y),\quad R(JX,Y)+R(X,JY)=0$$
for all vector fields on $M$. This implies

 \begin{lem} For all vector fields $Z,W\in T^{1,0}M$ it holds
$$R(Z,W)=R(\bar Z,\bar W)=0,\quad
 R(Z,\bar W)=\overline{R(\bar Z,W)},$$
 where the conjugation is taken with respect to
 $\u(T^\Co M)\subset\gl(T^\Co M)$.
 In particular, $$R_x(Z,\bar Z)\in i\g.$$
\end{lem}

It holds $$R(X,Y)=R(X,\bar Y)-R(Y,\bar X),$$ where on the left
hand side $R$ is the curvature tensor of the  connection in $TM$
and $X,Y$ are vector fields on $X$, and on the right hand side $R$
is the curvature tensor of the complexified connection and $X,Y$
are vector fields from $T^{1,0}M$ corresponding to the initial
$X,Y$.

Let now $\g\subset\u(1,n+1)_{\Co p}$ be an arbitrary subalgebra.
We may carry all the above notations to this case. Denote
$\Co^{1,n+1}$ simply by $V$, and as above we write
$$V^\Co=V\otimes \Co=V^{1,0}\oplus V^{0,1}=V\oplus\bar V.$$  For the Lie algebra
$\g^\Co=\g\otimes\Co$ consider the following space of algebraic
curvature tensors:
$$\R(\g^\Co)=\left.\left\{R:
\Lambda^2V^\Co\to\g^\Co\right|\begin{array}{c}R(X,Y)Z+R(Y,Z)X+R(Z,X)Y=0,\,\,\forall\,\,X,Y,Z\in
V^\Co\\ R(X,\bar Y)=\sigma(R(\bar X,Y)),\,\,\forall X,Y\in
V\end{array}\right\}.$$ We say that $\g\subset\u(1,n+1)_{\Co p}$
is a Berger algebra if $\g^\Co$ is generated by the images of the
elements $R: \Lambda^2V^\Co\to\g^\Co$ from the space $\R(\g^\Co)$.
From the Ambrose-Singer theorem and the above consideration it
follows that if $\g\subset\u(1,n+1)_{\Co p}$ is the holonomy
algebra of a Lorentz-K\"ahler manifold, then it is a Berger
algebra.

The proof of the following lemma is direct.
\begin{lem} A linear map $R:
\Lambda^2V^\Co\to\g^\Co$ belongs to $\R(\g^\Co)$ if and only if it
satisfies $$R(X,\bar Y)=\sigma(R(\bar X,Y)),\quad R(X,\bar
Y)Z=R(Z,\bar Y)X,\,\,\forall X,Y,Z\in V.$$
\end{lem}

For an element $\xi\in\g^\Co$ we denote by $\xi^{1,0}$ its
projection onto $\gl(V^{1,0})=\gl(V)$. The following lemma easily
follows from the previous one. We denote by $X,Y$ vectors from
$\Co^n$.

\begin{lem}\label{LemR} Each algebraic curvature tensor $R\in\R(\u(1,n+1)_{\Co
p }\otimes\Co)$ is uniquely given by the equalities
$$R^{1,0}(p,\bar q)=\left(\begin{array}{ccc} \alpha & N^t
&\beta\\0&0&0\\0&0&0\end{array}\right),\quad R^{1,0}(X,\bar q)=
\left(\begin{array}{ccc} g(N,X) & T(X)^t&g(\bar
K,X)\\0&P(X)&AX\\0&0&0\end{array}\right),$$

$$ R^{1,0}(X,\bar
Y)= \left(\begin{array}{ccc} 0 & g(P(X)\cdot,\bar Y) & g(AX,\bar
Y)\\0&R_0(X,\bar
Y)&\overline{P(Y)}^tX\\0&0&0\end{array}\right),\quad
R^{1,0}(q,\bar q)= \left(\begin{array}{ccc} \beta & \bar K^t
&c\\0&A&K\\0&0&\bar \beta\end{array}\right)$$
 for arbitrary
elements $\alpha,\beta,\in\Co$, $c\in\Real$, $N\in
\overline{\Co^n}$, $K\in\Co^n$, $T\in\odot^2\overline{\Co^n}$,
$R_0\in\R(\u(n)^\Co)$, $A\in \gl(n,\Co)$ and
$P\in\gl(n,\Co)^{(1)}=\odot^2\overline{\Co^n}\otimes\Co^n$.
\end{lem}

Now, for an arbitrary subalgebra $\g\subset\u(1,n+1)_{\Co p}$ it
holds $$\R(\g^\Co)=\R(\u(1,n+1)_{\Co p
}\otimes\Co)\cap(\Lambda^2(V^\Co)^*\otimes\g^\Co).$$ Recall that
we assign to each weakly irreducible subalgebra
$\g\subset\u(1,n+1)_{\Co p}$ a real vector subspace $L=\Co^m\oplus
L_0\subset \Co^n$ such that $L_0\subset\Co^{n-m}$ is a real  form.

\begin{theorem}\label{ThBerger}
Let $n\geq 1$. If $\g\subset\u(1,n+1)_{\Co p}$ is a weakly
irreducible Berger subalgebra, then it is one of the following
algebras:
\begin{itemize}

\item $\g^\k=\k\zr(L\zr i\Real)$\\
$=\left.\left\{\left(\begin{array}{cccc}
a_1+ia_2&-\bar Z^t&-\bar X^t &ic\\
0&A&0&Z\\
0&0& a_2(iE_{n-m}+\theta)&X\\
0&0&0&-a_1+ia_2\end{array}\right)\right|\begin{array}{c}a_1+a_2(i+i\id_{\Co^{n-m}}+\theta)+A\in\k,\\
Z\in\Co^{m},\,X\in L_0,\\ c\in\Real\end{array}\right\},$

where
$$\k\subset\Real\oplus\Real(i,i\id_{\Co^{n-m}}+\theta)\oplus\u(m)$$ is an
arbitrary subalgebra;

\item $\g^{\k,\psi}=\k\oplus\{(0,\psi(X),X,0)|X\in \Co^{m-r}\oplus L_0\}\zr(\Co^r\zr i\Real)$
$$=\left.\left\{\left(\begin{array}{cccc}
0&-\bar Z^t&-\bar X^t &ic\\
0&A+\psi(X)&0&Z\\
0&0& 0&X\\
0&0&0&0\end{array}\right)\right|\begin{array}{c}A\in\k,\\
Z\in\Co^{r},\,X\in \Co^{m-r}\oplus L_0,\\
c\in\Real\end{array}\right\},$$ where $1\leq r\leq m$, we consider
the decomposition
$$\Co^m=\Co^{r}\oplus\Co^{m-r},$$ $\k\subset\u(r)$ is an arbitrary
subalgebra,
$$\psi:\Co^{m-r}\oplus L_0\to \u(r)$$ is a non-zero linear map such
that the image $\psi(\Co^{m-r}\oplus L_0)\subset\u(r)$ is
commutative, it commutes with $\k_0$ and has trivial intersection
with $\k_0$.
\end{itemize}
\end{theorem}

\subsection{Proof of Theorem \ref{ThBerger}}

We consider the Lie algebras obtained in Section \ref{secW-Ir} and
verify, which of these Lie algebras  are Berger algebras.

\begin{lem}  If $\g\subset\u(1,n+1)_{\Co p}$ is a
weakly irreducible Berger subalgebra contained in
$$\Real\oplus\Real(i+i\id_{\Co^{n-m}})\oplus\u(m)\oplus\z_{\u(n-m)}\tau\zr (L\zr
i\Real),$$ then $\g$ is contained in
$$\Real\oplus\Real(i+i\id_{\Co^{n-m}}+\theta)\oplus\u(m)\zr (L\zr
i\Real).$$
\end{lem}

{\bf Proof.} Let $R\in\R(\g^\Co)$ be as in Lemma \ref{LemR}. It is
clear that both $R_0$ and $P$ take their values in $\u(m)$.
Consider the tensor $A$, which should be of the form
$A=iA_1+iA_2+a_2\id_{\Co^{n-m}}$, where $A_1\in \u(m)$, $A_2\in
\z_{\u(n-m)}\tau$, and $a_2\in\Real$. Note that the
complexification of the vector subspace $L\subset\g$ has the form
\begin{equation}\label{compL}L\otimes\Co=\left.\left\{
\left(\begin{array}{cccc} 0 & \bar W^t &
-\overline{\tau(X)}^t&0\\0&0&0&Z\\0&0&0&X\\0&0&0&0\end{array}\right)\right|W,Z\in\Co^m,\,\,X\in\Co^{n-m}\right\}.\end{equation}
Let $X\in iL_0$, then this and  the condition $R(X,\bar
q)\in\g^\Co$ imply
$$T(Y,X)=-g(Y,\overline{\tau(iA_2X+a_2X)})$$
for all $Y\in\Co^n$. Since $X\in iL_0$ and $A_2$ preserves $iL$,
it holds $\tau(iA_2X+a_2X)=iA_2X-a_2X$. Using \eqref{g=4h} and the
fact that $T$ is symmetric, we obtain
$$h(Y,iA_2X-a_2X)=h(X,iA_2Y-a_2Y)$$
for all $Y\in iL_0$. This implies
$${\rm Re}h(X,A_2 Y)=a_2{\rm Im}
h(X,Y)$$ for all $X,Y\in L_0$, i.e. $A_2=-a_2\theta$. This proves
the lemma. $\Box$

\begin{lem}  If  $\g\subset\u(1,n+1)_{\Co p}$ is a
weakly irreducible Berger subalgebra, then $i\Real\subset\g$.
\end{lem}

{\it Proof.} Let $\g\subset\u(1,n+1)_{\Co p}$ be a weakly
irreducible Berger subalgebra and suppose that
$i\Real\not\subset\g$. Then it is clear that $m=0$, i.e. $L=L_0$.
By the previous lemma, $\pr_{\Co\oplus\u(n)}\g\subset
\Real(i+i\id_{\Co^{n-m}}+\theta)$. Let $R\in\R(\g^\Co)$ be as in
Lemma \ref{LemR}. From the condition that $R$ takes values in
$\g^\Co$ it follows immediately that $R=0$. \qed

\begin{lem} If  $\g\subset\u(1,n+1)_{\Co p}$ is a
weakly irreducible Berger subalgebra of the form
$(\f_0\oplus\{\psi(X)+X|X\in \Co^{m-r}\oplus L_0\})\zr (\Co^r\zr
i\Real)$, then $\f_0\oplus\psi(\Co^{m-r}\oplus L_0)\subset\u(r)$.
\end{lem}

{\it Proof.} Suppose that  $\pr_{\Co\oplus \u(n)}\g$ contains an
element of the form $i+i\id_{\Co^{n-m}}+A$, where $A\in\u(r)$.
Since $\g$ is a Berger algebra, there exists $R\in\R(\g^\Co)$ such
that
$$R^{1,0}(q,\bar q)=\left(\begin{array}{ccc} 1 & *
&*\\0&\id_{\Co^{n-m}}-iA&*\\0&0&1\end{array}\right).$$
 Let
$X\in\Co^{m-r}\oplus L_0$. It holds
$$R^{1,0}(X,\bar q)=\left(\begin{array}{cccc} 0 & * &
\pr_{\Co^{m-r}\oplus
L_0}T(X)^t&*\\0&0&0&*\\0&0&0&X\\0&0&0&0\end{array}\right).$$ This
implies that $X-\pr_{\Co^{m-r}\oplus L_0}\overline{T(X)}\in\g$.
Considering $iX$, we get $iX+i\pr_{\Co^{m-r}\oplus
L_0}\overline{T(X)}\in\g$. The vectors of the last two types span
$\Co^{m-r}\oplus L_0$. We conclude that $\psi=0$. This gives a
contradiction. \qed

The lemmas show that if $\g\subset\u(1,n+1)_{\Co p}$ is a weakly
irreducible Berger subalgebra, then $\g$ must be one of the
algebras from the statement of the theorem. Conversely it is easy
to see that all algebras from the statement of the theorem are
weakly irreducible Berger algebras. \qed

\section{Lorentz-K\"ahler metrics}\label{secLKmetrics}

\subsection{Complex Walker metrics}\label{secWalker}
Suppose that $(M,g)$ is a Lorentzian manifold such that its
holonomy group preserves an isotropic line, then $(M,g)$ admits a
parallel distribution $\ell$ of isotropic lines. Locally there
exist so called Walker coordinates $v,x^1,...,x^n,u$ such that the
metric $g$ has the form $$g=2dvdu+h+2Adu+H (d u)^2,$$ where
$h=h_{ij}(x^1,...,x^n,u)d x^id x^j$ is an $u$-dependent family of
Riemannian metrics,  $A=A_i(x^1, \ldots, x^n,u)d x^i$ is an
$u$-dependent family of one-forms, and $H=H(v,x^1,...,x^n,u)$ is a
local function on $M$ \cite{WalkerBook}.

Let us show that Walker coordinates may be generalized to the case
of a Lorentz-K\"ahler manifold. Suppose that $(M,h)$ is a
Lorentzian-K\"ahler manifold such that its holonomy group preserves
a complex isotropic line in the tangent space. Then $(M,h)$ admits a
parallel distribution $\ell\subset TM$ of complex isotropic lines.
It is clear that the corresponding line bundles $\ell\subset
T^{1,0}M$ and $\bar\ell\subset T^{0,1}M$ are parallel. Moreover, the
distribution $\ell^\bot\subset TM$ of complex codimension 1 is
parallel, and the subbundles $\ell^\bot\subset T^{1,0}M$ and
$\overline{\ell^\bot}\subset T^{0,1}M$ are parallel as well. It is
clear that all parallel distributions are involutive. From the
holomorphic version of the Frobenius Theorem it follows that locally
on $M$ there exist  holomorphic coordinates $v,z^1,\dots,z^n,u$ such
that the vector field $\partial_v$ generates the subbundle
$\ell\subset T^{1,0}M$, and the vector fields
$\partial_v,\partial_{z^1},\dots,\partial_{z^n}$ generate the
subbundle $\ell^\bot\subset T^{1,0}M$. Consequently, the metric $h$
can be written in the following way:
$$h=h_{\bar uv}d\bar udv+h_{\bar v u}d\bar vdu+h_{\bar jk}d\bar
z^jdz^k+h_{\bar uk}d\bar udz^k+h_{\bar k u}d\bar z^kdu+h_{\bar
uu}d\bar u du.$$ Since the metric $h$ is pseudo-K\"ahler, it holds
\begin{equation}\label{habhba}
\partial_a h_{\bar b c}=\partial_c h_{\bar b a},
 \quad \partial_{\bar a} h_{\bar b c}=\partial_{\bar b} h_{\bar a c}.\end{equation}
 Hence the coefficients of the metric depend on the
coordinates in the following way:
$$h_{\bar  v u}=h_{\bar  v u}(\bar v,\bar z^l,u,\bar u),\quad
h_{\bar  j k}=h_{\bar jk}(z^l,\bar z^l,u,\bar u),$$
$$h_{\bar   k u}=h_{\bar  k u}(\bar v,\bar z^l,u,\bar u),\quad
h_{\bar  uu}=h_{\bar  u u}(v,\bar v,z^l,\bar z^l,u,\bar u)=h^1_{\bar
u u}(v,z^l,u,\bar u)+h^2_{\bar  u u}(\bar v,\bar z^l,u,\bar
u)+h^3_{\bar  u u}(z^l,\bar z^l,u,\bar u).$$ Thus any such metric is
given by a potential $f$ that satisfies the equations
$$\partial_v\partial_{\bar v} f=\partial_v\partial_{\bar z^k} f=0.$$
It is easy the check that the inverse matrix, defined by the
equality
$$h^{\bar a b}h_{\bar ac}=\delta^b_c,$$
is given by $$h^{\bar vv}=\frac{1}{|h_{\bar uv}|^2}\left(h_{\bar
lu}h^{\bar lj}h_{\bar uj}-h_{\bar uu}\right),\quad h^{\bar
vu}=\frac{1}{h_{\bar v u}},\quad h^{\bar jk}=\tilde h^{\bar
jk},\quad h^{\bar vk}=-\frac{1}{h_{\bar v u}}h_{\bar j u}h^{\bar j
k},\quad h^{\bar uk}=h^{\bar uu}=0,$$ where $\tilde h^{\bar jk}$
is the inverse matrix to $h_{\bar jk}$.

For the obtained metric and $a\neq v$ it holds
$$\Gamma^a_{vb}=h^{\bar ca}\partial_b h_{\bar c v}=h^{\bar ua}\partial_b h_{\bar u v}=0,$$
i.e. the vector field $\partial_v$ is isotropic and recurrent,
consequently it defines a local parallel distribution of complex
isotropic lines $\ell\subset T^{1,0}M$.

We will consider the frame $p,e_1,\dots, e_n,q,\bar p,\bar
e_1,\dots, \bar e_n,\bar q$, where
$$p=\frac{1}{h_{\bar vu}}\partial_v,\quad e_j=C^k_j\left(\partial_{z^k}
-\frac{h_{\bar u k}}{h_{\bar u v}}\partial_v\right),\quad
q=\partial_u-\frac{h_{\bar uu}}{h_{\bar uv}}\partial_v.$$ Here
$C^k_j$ is a matrix such that $\bar C^k_j h_{\bar k l}
C^l_s=\delta_{jl}$. It is clear that such matrix exists.

\subsection{Non-existence result}\label{secNon-ex}

Here we will show that some of the Berger algebras obtained above
cannot appear as the holonomy algebras of Lorentz-K\"ahler
metrics.

\begin{theorem}\label{ThnonEx} Let $\g^\k$ be the Berger algebra
as in Theorem~\ref{ThBerger}. If $m<n$ and $\g^\k$ is the holonomy
algebra of a Lorentz-K\"ahler manifold, then $\k\subset
\Real(i+i\id_{\Co^{n-m}}) \oplus\u(m)$. Consequently, if
$\k\not\subset\u(m)$, then $\theta=0$, and
$L_0=\Real^{n-m}$.\end{theorem}

{\bf Proof.}  Consider the Berger algebra $\g^\k$ as in
Theorem~\ref{ThBerger}. Suppose that $m<n$  and suppose that
$\g^\k\otimes\Co$ is the complexification of the holonomy algebra
of a Lorentz-K\"ahler manifold $(M,h)$ at a point $x\in M$.
Suppose that $\k\otimes\Co$ contains an element
$$\xi=ia_1+iA+a_2(-\id+i\theta),$$ where $A\in\u(m)$. Suppose that
$a_1\neq 0$, or $a_2\theta\neq 0$. Consider the coordinates and a
local frame as in the previous section. Note that the holonomy
algebra of the induced connection on the bundle $\ell^\bot/\ell$
coincides with $\pr_{\u(n)}\k$; this subbundle and its connection
may be identified with the subbundle of $T^{1,0}M$ generated by
the vector fields $e_1,\dots,e_n$ and the connection $\nabla$
restricted and then projected to that subbundle. Let $\lambda$ be
one of the non-zero eigenvalues of $i\theta$. Let $E_{\lambda
x}\subset \spa_{\Co}\{e_{m+1},\dots,e_n\}$ be the eigenspace
corresponding to the eigenvalue $-a_2(1-\lambda)$ of the operator
$\xi$ restricted to $\spa_{\Co}\{e_{m+1},\dots,e_n\}$. The
corresponding subspace in $\ell_x^\bot/\ell_x$ is clearly
holonomy-invariant, and we obtain a parallel subbundle
$E_\lambda\subset \ell^\bot/\ell$. Denote by the same symbol
$E_\lambda$ the corresponding subbundle of $T^{1,0}M$, which is
parallel modulo $p$. The holonomy algebra $\g^\k\otimes\Co$
preserves the space of tensors of the form
$$\bar p_x\wedge Z-p_x\wedge \overline{\kappa(Z)},$$
where $Z\in E_{\lambda x}$ and $\kappa(Z)$ may be found from the
condition
$$h(\overline{\tau(X)},Z)+h(\overline{\kappa(Z)},X)=0,\quad
\forall X\in\spa_{\Co}\{e_{1x},\dots, e_{nx}\}.$$ The element $\xi$
acts in this subspace as the multiplication by $a_1i+a_2\lambda$,
and its orthogonal complement in the holonomy algebra acts in this
space trivially. Consequently we get a parallel subbundle
$F_\lambda$ of the bundle of bivectors. Let $y\in M$, and $\gamma$
be any curve from $x$ to $y$. Then $F_{\lambda y}=\tau_\gamma
F_{\lambda x}$, and the element
$\tau_\gamma\circ\xi\circ\tau^{-1}_\gamma$ of the holonomy algebra
at the point $y$ acts on $F_{\lambda y}$ as the multiplication by
$a_1i+a_2\lambda$. The same element acts on $\Co p_y\wedge \bar p_y$
as the multiplication by $2a_1i$. Consequently, $F_\lambda$ has
trivial projection to the bundle $<p\wedge \bar p>$ generated by
$p\wedge \bar p$, and $F_\lambda$ consists of sections of the form
$$\bar p\wedge Z-p\wedge \bar W,$$ where $Z\in\Gamma(E_\gamma)$
and $W$ is uniquely defined by $Z$. For each such section and each
vector field $V$ from $T^\Co M$ must hold
$$\pr_{<p\wedge \bar p>}\nabla_V(\bar p\wedge Z-p\wedge \bar
W)=0.$$ Hence, $$\bar p\wedge\pr_{<p>}\nabla_VZ-p\wedge\pr_{<\bar
p>}\nabla_V\bar W=0,$$ i.e.
$$h(\nabla_VZ,\bar q)+h(\nabla_V\bar W,q)=0.$$ Let
$V=\partial_{z^a}$. Since $\bar q=\partial_{\bar u}-\frac{h_{\bar
uu}}{h_{\bar vu}}\partial_{\bar v}$, it holds
$$h(\nabla_{\partial_{z^a}} Z,\bar q)=-h(Z,\nabla_{\partial_{z^a}}\bar
q)=0.$$ Let $V=\partial_{\bar z^a}$. By the similar argument it
holds $h(\nabla_{\partial{\bar z^a}}\bar W,q)=0$. Consequently,
$$h(\nabla_{\partial{\bar z^a}}Z,\bar q)=0.$$
This shows that the subbundle $E_\lambda\subset T^{1,0}M$ is
parallel. The corresponding distribution $E_\lambda\subset TM$ is
also parallel, and it is non-degenerate. This gives a
contradiction, since the holonomy algebra is weakly irreducible.
Thus, $a_1=0$, and $a_2\theta =0$. This proves the theorem. \qed

\subsection{Construction of the metrics with all possible holonomy
algebras}\label{secConstr}

{\bf 1)} First consider the complex dimension $2$.
 Metrics with the holonomy algebra $\g_0$ may be found in
\cite{BB-I97,G-T01}.

Let $a,b\in \Co$ be some numbers. We assume that if $a=0$, then
$b=0$. If $a=b=0$, then we set
$$f_\Co=\bar v u+\bar u v;$$
if $a\neq 0, b=0$, we set
$$f_\Co(v,\bar v,u,\bar u)=\left\{\begin{array}{ll} 2{\rm Re}\left(-\frac{v}{iau}\left(e^{-ia|u|^2}-1\right)
\right), & \text{ if } u\neq 0;\\
0, & \text{ if } u=0.\\
\end{array}\right.$$
We claim that the function $f_\Co$ is real analytic. Indeed, the
function $$F(z)=e^{-iaz}-1$$ is holomorphic and it vanishes at the
point $z=0$. This implies that there exists an holomorphic function
$H(z)$ such that
$$F(z)=zH(z).$$ Consequently,
$$\frac{1}{u}\left(e^{-ia|u|^2}-1\right)=\bar u H(|u|^2).$$
This shows that the function $f_\Co$ is real analytic. Note that it
holds
$$\partial_{\bar u}\partial_vf_\Co=e^{-ia|u|^2}.$$

Next, if $a,b\neq 0$, then let
$$f_\Co(v,\bar v,u,\bar u)=\left\{\begin{array}{ll} 2{\rm Re}\left(
e^{\frac{ia^2}{b}}\frac{\sqrt{\pi}v}{\sqrt{ib}u}\left({\rm
erf}\left(\frac{ia}{\sqrt{ib}}+\frac{\sqrt{ib}|u|^2}{2}\right)- {\rm
erf}\left(\frac{ia}{\sqrt{ib}} \right)\right)\right), & \text{ if }
u\neq 0;
\\ 0, & \text{ if } u=0.\\
\end{array}\right.$$
Recall that the Error function ${\rm erf}(z)$ satisfies
$$\frac{d}{dz}{\rm erf} (z)=\frac{2}{\sqrt{\pi}}e^{-z^2}.$$ The
function $f_\Co$ is again real analytic and it satisfies
$$\partial_{\bar u}\partial_vf_\Co=e^{-ia|u|^2-\frac{ib}{4}|u|^4}.$$
The function $f_\Co$ will be a part of the potential $f$ that we are
going to construct; this function  is chosen in such a way that the
Lorentz-K\"ahler matric for the potential $f$ satisfies the
condition $$h_{\bar uv}=e^{-ia|u|^2-\frac{ib}{4}|u|^4}.$$ We could
define the metric $h$ directly without considering its potential,
then we could avoid the usage of the Error function and functions
with $u$ in denominator, but then a priori the metric would be only
pseudo-Hermitian and we would have to check the condition
\eqref{habhba} in order to show that the metric is pseudo-K\"ahler.

\medskip

Let $\g$ denote the holonomy algebra of a metric given by some
potential $f$.  Similarly as in the proof of Theorem
\ref{ThConstruct} the following may be shown:

if we chose $f=f_{\Co}$ with $a=i$, $b=1$,  then $\g=\g_1$;

if we choose $f=f_{\Co}$ with $a=\gamma\neq 0$, $b=0$, then
$\g=\g_3^\gamma$;

if we choose $f=\bar uv+\bar vu+|u|^4$, then $\g=\g^0_3$;

finally, the holonomy algebra of the metric $h=e^{\bar u v}d\bar
udv+e^{\bar v u}d\bar vdu$ coincides with $\g_2$.

{\bf 2)} Let $n\geq 1$.  Denote by $v,z^1,...,z^n,u$ the
coordinates on $\Co^{n+2}$. Let $\k\subset\Co\oplus\u(n)$ be a
subalgebra. Let us choose the following elements spanning $\k$:
$$a+A_1,\quad b+A_2,\quad A_3,...,A_N,$$
where $a,b\in\Co$, $A_1,...,A_N\in\u(n)$. We may assume the
following: if $a=0$, then $b=0$; if $a\neq 0$ (resp. $b\neq 0$),
then $A_1$ (resp. $A_2$) belongs to the center of $\pr_{\u(n)}\k$.

Consider the following Hermitian matrices:
$$B_\alpha=-\frac{1}{(\alpha!)^2}iA_\alpha,\,\alpha=1,...,N$$ and
$$G=\sum_{\alpha=1}^NB_\alpha |u|^{2\alpha}.$$
Define the function $$f_{\u(n)}=\bar Z^Te^GZ,$$ where $e^G$ is the
matrix exponential and  $Z=\left(\begin{array}{c}z^1\\\vdots \\
z^n\end{array}\right)$.


We may suppose that the basis $e_1,\dots,e_n$ is chosen in such a
way that, for some $0 \leq n_0\leq m$, $A_1$ restricted to
$\Co^{n_0}=\spa_{\Co}\{e_1,\dots,e_{n_0}\}$ is an isomorphism of
$\Co^{n_0}$, and it annihilates
$\spa_{\Co}\{e_{n_0},\dots,e_{n}\}$.
 Define the following function:
$$f_{\Co^m}=\frac{1}{2}{\rm Re}\left(i\bar
u^2\sum_{k=n_0+1}^m(z^k)^2\right).$$

Define the function
$$f_1=f_\Co+f_{\u(n)}+f_{\Co^n}.$$

\medskip

Consider the Lie algebra $\g^{\k,J,L}$. We assume that $a=i$,
$b=0$,
$$A_1=\tilde A_1+i\id_{\Co^{n-m}},$$
where $\tilde A_1\in u(m)$, and $A_2,\dots,A_N\in \u(m)$. Let
$n_0$ be the number defined by $\tilde A_1$ in the same way as it
was defined above by $A_1$. Let
$$f_{\Real^{n-m}}(z^l,\bar z^l,u,\bar u)=\left\{\begin{array}{ll}- {\rm Re}\left(
\sum_{j=m+1}^n(\bar z^j)^2\frac{1}{\bar
u^2}\left(1-e^{|u|^2}+|u|^2e^{|u|^2}\right)\right), & \text{ if }
u\neq 0;
\\ 0, & \text{ if } u=0.\\
\end{array}\right.$$
Note that the function $$F(z)=1-e^z+ze^z$$ is holomorphic and it
satisfies $F(0)=\frac{d}{dz}F(0)=0$. This implies that there exists
a holomorphic function $G(z)$ such that $F(z)=z^2G(z)$. Consequently
the function $f_{\Real^{n-m}}$ is real analytic. It holds
$$\partial_{\bar z^j}\partial_u f_{\Real^{n-m}}=-\bar z^j
ue^{|u|^2},\quad j=m+1,\dots,n.$$ Define the function
$$f_2=f_\Co+f_{\u(n)}+f_{\Co^m}+f_{\Real^{n-m}}.$$

\medskip

Consider the Lie algebra $\g^{\k,L}$. We assume that $a=b=0$,
$A_1,\dots,A_N\in \u(m)$. Using the real form
$L_0\subset\Co^{n-m}$, we define the matrix
$$B=(B_{kj})_{k,j=m+1}^n={\rm
diag}\left(\frac{\sqrt{2}}{2}\left(\begin{matrix}\sqrt{1-\lambda_1}&\sqrt{1+\lambda_1}\\
-i\sqrt{1+\lambda_1}&-i\sqrt{1-\lambda_1}\end{matrix}\right),\dots,\frac{\sqrt{2}}{2}\left(\begin{matrix}\sqrt{1-\lambda_s}&\sqrt{1+\lambda_s}\\
-i\sqrt{1+\lambda_s}&-i\sqrt{1-\lambda_s}\end{matrix}\right),1,\dots,1\right).$$
Let $$f_{L_0}=-2{\rm
Re}\left(\sum_{j=m+1}^n\sum_{\alpha=1}^{n-m}\frac{iB_{j\,m+\alpha}\bar
z^j}{((N+\alpha)!)^2}|u|^{2(N+\alpha)}
\frac{u}{N+\alpha+1}\right),$$ and
$$f_3=f_\Co+f_{\u(n)}+f_{\Co^m}+f_{L_0}.$$

\medskip

Let us turn to the algebra $\g^{\k_0,\psi}$. Let
$\k=\k_0\oplus\psi(\Co^{m-r}\oplus L_0)$.
 Consider the unitary
matrices
$$A_1=\psi(e_{r+1}),...,A_{m-r}=\psi({e_m}),\quad
A_{m-r+1}=\psi(ie_{r+1}),...,A_{2m-2r}=\psi(i{e_m}),$$
$$A_{2m-2r+1}=\psi(f_{m+1}),...,A_{n+m-2r}=\psi(f_n),$$
where, as above, $f_{m+1},\dots,f_n$ is a basis of the real vector
space $L_0$. Let $A_{n+m-2r+2},...,A_N$ be a basis of $\k_0$. Let
$$D=\left(\begin{array}{ccc}i
E_{m-r}&-E_{m-r}&0\\0&0&iB\end{array}\right).$$ Using this matrix,
we define the function
$$f_\psi=-{\rm
Re}\left(\sum_{j=r+1}^n\sum_{\alpha=1}^{n+m-2r}\frac{D_{j\,\alpha}(\bar
z^j)^2}{((\alpha+2)!)^2}|u|^{2(\alpha+2)}
\frac{u}{\alpha+3}\right).$$ Finally consider the function
$$f_4=f_\Co+f_{\u(n)}+f_{\Co^r}+f_{\psi}.$$

\begin{theorem}\label{ThConstruct} The holonomy algebras of
the Lorentz-K\"ahler metrics with the potentials $f_1$, $f_2$,
$f_3$, $f_4$ coincide respectively with the Lie algebras $\g^\k$,
$\g^{\k,J,L}$, $\g^{\k,L}$, $\g^{\k_0,\psi}$.
\end{theorem}

\subsection{Proof of Theorem \ref{ThConstruct}}\label{secProofConstr}

Let $(M,h)$ be a Lorentz-K\"ahler manifold, let $\g\subset\u(T_xM)$
be the holonomy algebra of the Levi-Civita connection on $(M,h)$.
Recall that the complexification $\g\otimes\Co\subset\gl(T^\Co_xM)$
coincides with holonomy algebra of the induced complex connection on
$T^\Co M$, and it is determined by the projection of $\g\otimes\Co$
to $\gl(T^{1,0}_xM)$. If the pseudo-Riemannian metric corresponding
to $h$ is analytic, then the holonomy algebra $\g\otimes\Co$ is
generated by the following elements:
\begin{equation}\label{nablasR}\nabla_{Z_r}\cdots\nabla_{Z_1}R_x(X,\bar Y),\quad r\geq 0,\,\,
X,Y\in T^{1,0}_xM,\,\,Z_1,\dots, Z_r\in T_x^\Co M.\end{equation} Let
$z^a$ be local complex coordinates on $M$. We denote by
$R^{1,0}(\partial_{z^a},\partial_{\bar z^b})$ the  matrix of the
field of the endomorphisms $R(\partial_{z^a},\partial_{\bar z^b})$
restricted to $T^{1,0}M$. Recall that the Christoffel symbols
$\Gamma^A_{BC}$ are given by the equality $$\nabla_{
\partial_{z^B}}{\partial_{z^A}}=\Gamma^A_{CB}{\partial_{z^C}},$$ where we
assume that $A,B,C$ take all values $a$ and $\bar a$. The possibly
non-zero symbols are $\Gamma^a_{bc}$ and $\Gamma^{\bar a}_{\bar
b\bar c}=\overline{\Gamma^a_{bc}}$. It holds
$$\Gamma^a_{bc}=h^{\bar da}\partial_ch_{\bar d b}.$$
The components of the curvature tensor are defined by the equality
$$R(\partial_{z^C},\partial_{z^D})\partial_{z^B}=R^A_{BCD}\partial_{z^A}.$$
Only the following coefficients may be  non-zero:
$$R^a_{bc\bar d}=-R^a_{b\bar d c},\quad -R^{\bar a}_{\bar b d \bar
c}= R^{\bar a}_{\bar b  \bar c d}=\overline{R^a_{bc\bar d}},$$ and
it holds
$$R^a_{bc\bar d}=-\partial_{\bar z^d}\Gamma^a_{bc}.$$ Let for each
fixed $c$,  $\Gamma_c$ denote the matrix $(\Gamma^a_{bc})$. Let
$\xi$ be the matrix of a one of the fields
$$\nabla_{\partial_{z^{A_r}}}\cdots\nabla_{\partial_{z^{A_1}}}R^{1,0}(\partial_{z^{a}},
\partial_{\bar z^{b}}).$$ It holds
\begin{equation}\label{nablaxi}\nabla_{\partial_{z^c}}\xi=\partial_{z^c}\xi+[\Gamma_c,\xi],\quad
\nabla_{\partial_{\bar z^c}}\xi=\partial_{\bar
z^c}\xi.\end{equation}

Let now $h$ be the metric defined by a potential $f$ from the
theorem. The coefficients of the metric $h$ are analytic as
complex-valued functions of real variables, i.e. the holonomy
algebra $\g\otimes\Co$ of the metric is generated by the elements
\eqref{nablasR}.

 It holds $$h_{\bar uv}=e^{-ia|u|^2-\frac{ib}{4}|u|^4},\quad
h_{\bar j k}=\left(e^G\right)_{\bar j k}.$$ We will consider the
frame $p,e_1,\dots, e_n,q,\bar p,\bar e_1,\dots, \bar e_n,\bar q$,
where
$$p=\frac{1}{h_{\bar uv}}\partial_v,\quad e_j=\left(e^{-\frac{1}{2}\bar G} \right)_{kj}\left(\partial_{z^k}
-\frac{h_{\bar u k}}{h_{\bar u v}}\partial_v\right),\quad
q=\partial_u-\frac{h_{\bar uu}}{h_{\bar uv}}\partial_v.$$ Using
this frame, the curvature tensor of the metric $h$ can be
expressed through some tensor fields in the same way as in Section
\ref{secBerger} (we will use the same notations).  Note that
$$\Gamma^v_{vv}=\Gamma^v_{vk}=\Gamma^{l}_{jk}=0.$$
Consequently, $\alpha=0$, $N=0$, $P=0$, and $R_0=0$.

Next, $$\Gamma^k_{uu}=h^{\bar lk}h_{\bar
vu}\partial_u\left(\frac{h_{\bar lu}}{h_{\bar vu}}\right),$$ and
\begin{multline*}K=\pr_{<e_1,\dots, e_n>}R(q,\bar
q)q=\pr_{<e_1,\dots, e_n>}R(\partial_u,\partial_{\bar
u})\partial_u\\ = \pr_{<e_1,\dots, e_n>}R^k_{uu\bar
u}\partial_{z^k}=\sum_{j=1}^n\left(e^{\frac{1}{2}\bar
G}\right)_{\bar jk}R^k_{uu\bar u}e_j.\end{multline*} Consequently,
$$K=-\sum_{k,j=1}^n\left(e^{\frac{1}{2}\bar
G}\right)_{\bar jk}\partial_{\bar
u}\left(\left(e^{-G}\right)_{\bar l k}h_{\bar
vu}\partial_u\left(\frac{h_{\bar lu}}{h_{\bar
vu}}\right)\right)e_j.$$ Similarly,
$$T(e_j)=-\sum_{k,k_1,j_1=1}^n\left(e^{-\frac{1}{2}\bar
G}\right)_{\bar j_1 j}\left(e^{-\frac{1}{2}\bar G}\right)_{\bar k_1
k}h_{\bar uv}\partial_{\bar u}\left(\frac{1}{h_{\bar
uv}}\partial_{z^{j_1}}h_{\bar uk_1}\right)\bar e_k.$$

Let $\g\subset\u(1,n+1)$ be the holonomy algebra of the metric $h$
at the point $0$.

\begin{lem}\label{lammak} It holds $\pr_{\Co\oplus\u(n)}\g=\k$.
\end{lem}

{\it Proof.}

It is easy to check that
$$\Gamma^v_{vu}=-ia\bar u-\frac{ib}{2}|u|^2\bar u,$$
and
$$R^p_p(q,\bar q)=ia+ ib|u|^2.$$
Consequently, a non-trivial value at the point $0$ at the position
$(v,v)$ may have only the matrices of the  following covariant
derivative of $R$:
\begin{equation}\label{nablavv}\nabla_{\partial_u}\nabla_{\partial_{\bar u}}R^v_{v u\bar
u}(0)= \nabla_{\partial_{\bar u}}\nabla_{ \partial_u}R^v_{vu\bar
u}(0)=bi.\end{equation}

Consider the matrix $\tilde\Gamma_u=(\Gamma_{ju}^k)_{j,k=1}^n$.
Let $\tilde R(X,\bar Y)$ be the similar matrix obtained from $R$
and vector fields $X,Y\in T^{1,0}M$.
 Let
$X,Y\in\{p,e_1,\dots,e_n,q\}$, then  the last matrix $\tilde
R(X,\bar Y)$ is non-zero only for $X=Y=q$. Moreover, $$\tilde
R(q,\bar q)=\tilde R(\partial_u,\partial_{\bar u}).$$
  It holds
$$\tilde \Gamma_u=e^{-G}\partial_ue^G.$$ The formula for the
derivation of the matrix exponential reads \cite[Sec. 5, Th.
1.2]{Rossmann}:
$$e^{-G}\partial_ue^G=\sum_{k=0}^\infty\frac{(-1)^k}{(k+1)!}({\rm
ad}_G)^k\partial_u
G=\partial_uG-\frac{1}{2!}[G,\partial_uG]+\frac{1}{3!}[G,[G,\partial_uG]]-\cdots.$$
Next, $$\tilde R(q,\bar q)=- \partial_{\bar u}\tilde \Gamma_u.$$

\begin{lem} It holds $$\nabla^r_{\partial_u}\nabla^r_{\partial_{\bar u}}\tilde R(\partial_u,\partial_{\bar
u})(0)=iA_r,\quad r=0,\dots,N.$$
\end{lem}

{\it Proof.} From \eqref{nablaxi} it follows that
$$\nabla^r_{\partial_{\bar u}}\tilde
R(\partial_u,\partial_{\bar u})=-\partial^{r+1}_{\bar
u}\tilde\Gamma_u=\frac{1}{r!}iA_ru^r-\frac{1}{2!}\partial^{r+1}_{\bar
u}[G,\partial_uG]+\frac{1}{3!}\partial^{r+1}_{\bar
u}[G,[G,\partial_uG]]-\cdots.$$ Note that
 the covariant derivative
with respect to $\partial_{u}$ contains the Lie brackets with the
matrix $\tilde\Gamma_u$. On the other hand, $\tilde\Gamma_u$ is
dividable by $\bar u$, consequently
\begin{multline*}\nabla^r_{\partial_u}\nabla^r_{\partial_{\bar u}}\tilde
R(\partial_u,\partial_{\bar u})(0)=\partial^r_u\partial^r_{\bar
u}\tilde R(\partial_u,\partial_{\bar u})(0)\\=iA_r-
\frac{1}{2!}\partial^{r+1}_{\bar
u}\partial^r_u[G,\partial_uG](0)+\frac{1}{3!}\partial^{r+1}_{\bar
u}\partial^r_u[G,[G,\partial_uG]](0)-\cdots.\end{multline*} Note
that since all terms in $G$ contain the same powers of $u$ and
$\bar u$, it holds
$$\partial^{r+1}_{\bar
u}\partial^r_u[G,\partial_uG](0)=\sum_{\begin{smallmatrix}
r_1+r_2=r+1\\r_1\geq 0,\,r_2\geq
1\end{smallmatrix}}[\partial^{r_1}_{\bar
u}\partial^{r_1}_uG,\partial^{r_2}_{\bar
u}\partial^{r_2}_uG](0)=[G,\partial^{r+1}_{\bar
u}\partial^{r+1}_uG](0)=0.$$ Here we used the symmetry in $r_1$,
$r_2$ combined with the skew-symmetry of the Lie brackets, and the
fact that $G(0)=0$. By the same reason $\partial^{r+1}_{\bar
u}\partial^r_u[G,[G,\partial_uG]](0)$ as well as all terms
containing more Lie brackets with $G$, are zero.
 \qed

Note that, in particular,
$$\tilde
R(\partial_u,\partial_{\bar u})(0)=iA_1,\quad
\nabla_{\partial_u}\nabla_{\partial_{\bar u}}\tilde
R(\partial_u,\partial_{\bar u})(0)=\nabla_{\partial_{\bar u}}
\nabla_{\partial_u}\tilde R(\partial_u,\partial_{\bar u})(0)=iA_2.$$
The elements $A_1$ and $A_2$ are in the center of $\pr_{\u(n)}\k$,
consequently, by the construction, the other covariant derivatives
of $\tilde R$ at the point $0$ take values in the subalgebra of
$\pr_{\gl(n,\Co)}\k$ orthogonal to $A_1$ and $A_2$. Using this and
comparing the last equalities with \eqref{nablavv}, we conclude that
$\k\subset\pr_{\Co\oplus\u(n)}\g$. The inverse inclusion easily
follows from the construction and the above formulas. Lemma
\ref{lammak} is proved. \qed

Let now $f=f_1$.

\begin{lem}\label{lemCn} It holds $\Co^n\zr i\Real\subset\g$.
\end{lem}

{\it Proof.} We know already that each holonomy algebra contains
$i\Real$. It holds $$R^{1,0}(e_j,\bar
q)(0)=\left(\begin{array}{ccc} 0 &
T(e_j)^t&*\\0&0&A_1e_j\\0&0&0\end{array}\right)\in\g^{1,0}.$$
Next, $T(e_j)=0$ for $j=1,\dots,n_0$, and $T(e_j)=-i\bar e_j$ for
$j=n_0+1,\dots,n$.
 This shows that $\g$ contains $A_1e_j, iA_1 e_j\in\Co^n$ for $j=1,\dots, n_0$. From
 the definition of $n_0$ it follows that $\g$ contains
 $\Co^{n_0}$. For each $j>n_0$,  $\g$ contains $e_j$ and $ie_j$.
 We conclude that $\g$ contains $\Co^n\zr i\Real$. \qed

From Lemmas \ref{lammak} and \ref{lemCn} it follows that
$\g=\g^\k$.

Suppose that $f=f_2$.

\begin{lem}\label{lemCmR} It holds $\pr_{\Co^n\zr i\Real}\g=L\zr i\Real=\Co^m\oplus \Real^{n-m}\zr i\Real\subset\g$.
\end{lem}

{\it Proof.} Similarly as in Lemma \ref{lemCn} it can be shown
that $\Co^m\zr i\Real\subset\g$. Let $j\geq m+1$. It holds
$$R(e_j,\bar q)=
\left(\begin{array}{cccc} 0 & * & \bar
e_j^t&*\\0&0&0&*\\0&0&0&-e_j\\0&0&0&0\end{array}\right).$$
Consequently, $\g$ contains $\Real^{n-m}$. The covariant
derivatives of this tensor  give trivial projections on
$\Co^{n-m}$. Finally, $K=0$, i.e. $$R^{1,0}(q,\bar
q)=-\id_{T^{1,0}M}+c p\wedge \bar p$$ for some function $c$.
Hence, the covariant derivatives of $R^{1,0}(q,\bar q)$ give
trivial projections on $\Co^{n-m}$. \qed

From Lemmas \ref{lammak} and \ref{lemCmR} it follows that
$\g=\g^{\k,J,L}$.

Suppose that $f=f_3$. In this case $\k\subset\u(m)$. As above it
can be shown that $\Co^m\zr i\Real\subset\g$. It holds
$$\pr_{<e_{m+1},\dots ,e_n>}K=\sum_{\alpha=1}^{n-m}\frac{iB_{j\,m+\alpha}}{((N+\alpha-1)!)^2}|u|^{2(N+\alpha-1)}e_j.$$ Together with the
proof of Lemma \ref{lammak} this shows that
$\pr_{\u(n)}\g=\k\subset\g$. The covariant derivatives of $K$ will
give $L_0\subset\g$. Finally, $\pr_{<e_{m+1},\dots, e_n>}\circ
T=0.$ Thus, $\g=\g^{\k,L}$.

The metric with the potential $f_4$ may be considered in the same
way. We will get $\g=\g^{\k_0,\psi}$. The theorem is proved. \qed

\section{Example: the space of oriented lines in
$\Real^3$}\label{secLines}

The space of oriented lines in $\Real^3$ admits the following
Lorentz-K\"ahler metric \cite{G-K08}:

$$h=\frac{1}{(1+|u|^2)^2}\left(d\bar u dv+d\bar vdu+\frac{2i(\bar
vu+\bar uv)}{(1+|u|^2)^3}d\bar udu\right).$$

It is easy to show that it holds
$$R^{1,0}(p,\bar
q)(0)=\left(\begin{matrix}0&2\\0&0\end{matrix}\right),\quad
R^{1,0}(q,\bar
q)(0)=\left(\begin{matrix}2&*\\0&2\end{matrix}\right).$$ This
implies that the holonomy algebra of the metric $h$ coincides with
$\u(1,1)_{\Co p}$.

\section{Complex pp-waves}\label{secpp-w}

Here we give the following characterization of the complex
pp-waves.

\begin{theorem}\label{Thpp-wave}

Let $(M,h)$ be a Lorentz-K\"ahler manifold of complex dimension
$n+2\geq 2$ with a parallel isotropic vector field
$p\in\Gamma(TM)$. Then the following conditions are equivalent:

\begin{itemize}

\item[1)] The holonomy algebra $\g$ of $(M,h)$ is contained in
$\Co^n\zr i\Real\subset\u(1,n+1)_{\Co p}$.

\item[2)] The curvature tensor of the Levi-Civita connection
satisfies $R(p^\bot,p^\bot)=0$.

\item[3)] The curvature tensor of the  extension of the Levi-Civita
connection to $T^\Co M$ satisfies $R(p^\bot,\bar p^\bot)=0$.

\item[4)] Locally $M$ admits complex coordinates
$v,z^1,\dots,z^n,u$
such that
$$h=d\bar udv+d\bar vdu+\sum_{k=1}^nd\bar
z^kdz^k+h_{\bar uk}d\bar udz^k+h_{\bar k u}d\bar z^kdu+h_{\bar
uu}d\bar u du,$$ and
 the coefficients of the metric depend on the
coordinates in the following way:
$$h_{\bar   k u}=h_{\bar  k u}(\bar z^l,u,\bar u),\quad
h_{\bar  uu}=h_{\bar  u u}(z^l,\bar z^l,u,\bar u),\quad
\partial_{z^k}\partial_{\bar z^j}h_{\bar uu}=0.$$

\item[5)] Locally $M$ admits complex coordinates
$v,z^1,\dots,z^n,u$ such that the potential of $h$ is of the form
$$f=\bar u v+\bar vu+\sum_{k=1}^n|z^k|^2+{\rm Re}(\phi(z^j,u,\bar u)).$$
\end{itemize}
\end{theorem}

{\bf Proof.} From results of Section \ref{secBerger} it follows
that the first three conditions are equivalent.  It is obvious
that the last two conditions are equivalent. Using computations as
in Section \ref{secProofConstr}, it is not hard to show that the
fourth condition implies the third one. Let us show that the third
condition implies the fourth one.

Consider local coordinates as in Section \ref{secWalker}. There
exists a function $\alpha$ such that $p=\alpha \partial_v$. Since
$\nabla p=0$, it holds $\partial_{\bar z^a}\alpha=0$. Considering
the coordinate transformation $$\tilde v=F(v,z^k,u),\quad \tilde
z^k=z^k,\quad \tilde u=u,$$ where the function $F$ satisfies
$\partial_vF=\frac{1}{\alpha}$, we get $p=\partial_{\tilde v}$.
This allows us to assume that $p=\partial_v$. The condition
$\nabla\partial_v=0$ implies $$0=\Gamma^v_{va}=h^{\bar
uv}\partial_ah_{\bar uv}.$$ Consequently, $$\partial_ah_{\bar
uv}=0.$$ from this and Section \ref{secWalker} it follows that
$h_{\bar uv}$ is a function of $\bar u$. Introducing the new
coordinate $\tilde u$ such that $h_{\bar uv}d\bar u=d\bar{\tilde
u}$, we get $h_{\bar {\tilde u}v}=1$, i.e. we may assume that
$h_{\bar u v}=1$.

It is clear that the induced connection in the bundle $p^\bot/<p>$
(where $p$ is considered as a vector field in $T^{1,0}M$) is flat.
Consequently, there exist  vector fields $e_1,\dots,e_n,\bar
e_1,\dots,\bar e_n$ such that $h(\bar e_k, e_j)=\delta_{ij}$,
$$e_j=B_j^k\partial_{z^k}\quad \text{modulo } p;$$ moreover, these vector fields are parallel
 modulo $p$. From this condition it follows that
 $$\partial_{v}B_j^k=\partial_{\bar v}B_j^k=\partial_{\bar u}B_j^k=\partial_{\bar z^l}B_j^k=0.$$

Consider the family of K\"ahler metrics $h_0=h_{\bar kj}d\bar
z^kdz^j$ depending on the parameters $u,\bar u$.  The $u$-families
of vector fields $e_1,\dots,e_n,\bar e_1,\dots,\bar e_n$ are
parallel with respect to the connections defined by the families of
these metrics. Consequently,
$$[e_j,e_k]=[\bar e_j,e_{ k}]=[\bar e_{ j},\bar e_{ k}]=0,$$
and there exist coordinates $\tilde z^1,\dots \tilde z^n$ such
that $e_j=\partial_{\tilde z^j}$ and $$h_0=\sum_{j=1}^n d{\bar
z^j}dz^j.$$ Since $B^j_k$ are functions of the variables
$z^1,\dots z^n,u$, the coordinates $\tilde z^1,\dots, \tilde z^n$
are related to $z^1,\dots,  z^n$ by a holomorphic transformation
depending holomorphically on the parameter $u$. By this reason, we
may consider the new coordinates $v,\tilde z^1,\dots, \tilde
z^n,u$. With respect to these coordinates it holds $h_{\bar k
j}=0$. The equality \eqref{habhba} implies the proof of the
theorem. \qed

\section{Lorentz-K\"ahler symmetric spaces}\label{secSym}

Berger \cite{BerSym} classified indecomposable simply connected
pseudo-Riemannian symmetric spaces with simple groups of
isometries. Classification of simply connected Lorentz-K\"ahler
symmetric spaces is obtained in \cite{K-O09}. Here we give an
alternative formulation of this result in terms of the curvature
and holonomy.

Let $(M,g)$ be a simply connected pseudo-Riemannian symmetric
space. Let $H$ be the group of transvections and $G\subset H$ be
the stabilizer of a point  $x\in M$, then the holonomy group of
$(M,g)$ coincides with the isotropy representation of $G$. The
groups $H$ and $G$, and consequently the manifold $(M,g)$,  may me
reconstructed from the holonomy algebra $\g\subset\so(\m)$
($\m=T_xM$), of $(M,g)$ and the value $R_x$ of curvature tensor of
$(M,g)$ at the point $x$ in the following way, see e.g.
\cite{Al10}. Define on the vector space
$$\h=\g\oplus\m$$
 the structure of the Lie algebra in the following way:
$$[A,X]=AX,\quad [A,B]=[A,B]_\g,\quad [X,Y]=-R_x(X,Y),$$
where $A,B\in\g$ and $X,Y\in\m$. Then $H$ may be found as the
simply connected Lie group with the Lie algebra $\h$, and
$G\subset H$ is the connected Lie subgroup corresponding to the
subalgebra $\g\subset\h$. Note that it holds $\g=R_x(\m,\m)$.

Thus in order to classify indecomposable simply connected
Lorentz-K\"ahler symmetric spaces it is enough to find all weakly
irreducible holonomy algebras $\g\subset\u(\m)$ ($\m=\Co^{1,n+1}$)
admitting  $\g$-invariant surjective maps $R:\Lambda^2\m\to\g$
that satisfy the Bianchi identity. It is convenient to consider
the complex extension of $R$ and use the results of Section
\ref{secBerger}.

For an indecomposable symmetric space $(M,g)$ the following
conditions are equivalent \cite{Al10}: $\h$ is simple;
$\g\subset\so(\m)$ is totally reducible; the Ricci tensor of
$(M,g)$ is non-degenerate.

If $\g\subset\u(1,n+1)$ is irreducible, then the above equivalent
conditions imply that $g=\u(1,n+1)$, and we obtain the complex
de~Sitter and anti de~Sitter symmetric spaces:
$$\textrm{dS}^{n+2}(\Co)=\SU(1,n+2)/\Un(1,n+1),\quad
\textrm{AdS}^{n+2}(\Co)=\SU(2,n+1)/\Un(1,n+1).$$ The subalgebra
$\g_0\subset\u(1,1)$ cannot be the holonomy algebra of a symmetric
space since it is completely reducible and is contained in
$\su(1,1)$ (i.e. the corresponding space would have  to be
Ricci-flat).

Thus we may assume that $\g\subset\u(1,n+1)$ is weakly irreducible
and it is contained in $\u(1,n+1)_{\Co p}$. Using the
classification of holonomy algebras, results of Section
\ref{secBerger} and solving a simple exercise in linear algebra we
arrive to the following theorem (we give the non-zero values of
$R$ on the basis vectors).

\begin{theorem}
If $(M,g)$ is an indecomposable Lorentz-K\"ahler symmetric space
with the holonomy algebra $\g\subset\u(1,n+1)_{\Co p}$, then
$(M,g)$ is given by exactly one of the following pairs $(\g,R)$:

\begin{itemize}
\item[a)] $\g=\left.\left\{\left(\begin{array}{cc}
0&ic\\
0&0\end{array}\right)\right|c\in\Real\in\Co\right\}$ and
$R^{1,0}(q,\bar q)=\left(\begin{array}{cc}
0&1\\
0&0\end{array}\right)$;

\item[b)] $\g$ and $-R$, where $\g$ and $R$ are from a);

\item[c)] $\g=\left.\left\{\left(\begin{array}{cc}
a&0\\
0&-\bar a\end{array}\right)\right|a\in\Co\right\}$ and
$R^{1,0}(p,\bar q)=\left(\begin{array}{cc}
1&0\\
0&0\end{array}\right)$;

\item[d)] $n=1$, $\g=\left.\left\{\left(\begin{array}{ccc}
0&-x &ic\\
0&0&x\\
0&0&0\end{array}\right)\right|x, c\in\Real\right\}$ and

$R^{1,0}(e_1,\bar q)=\left(\begin{array}{ccc}
0&0 &-i\\
0&0&0\\
0&0&0\end{array}\right)$, $R^{1,0}(q,\bar
q)=\left(\begin{array}{ccc}
0&-i &0\\
0&0&i\\
0&0&0\end{array}\right)$;

\item[e)] $\g$ and $-R$, where $\g$ and $R$ are from d);

\item[f)] $n\geq 1$, $0\leq m\leq n$,
$\g=\left.\left\{\left(\begin{array}{cccc}
2ai&-\bar Z^t&-\bar X^t &ic\\
0&aiE_m&0&Z\\
0&0& 2aiE_{n-m}&X\\
0&0&0&2ai\end{array}\right)\right|\begin{array}{c}
Z\in\Co^{m},\,X\in \Real^{n-m},\\ a,c\in\Real\end{array}\right\},$

$R^{1,0}(p,\bar q)=2R^{1,0}(e_j,\bar
e_j)=\frac{1}{2}R^{1,0}(e_k,\bar e_k)=\left(\begin{array}{ccc}
0&0 &1\\
0&0&0\\
0&0&0\end{array}\right),$

$R^{1,0}(e_j,\bar q)=\frac{1}{2} \left(\begin{array}{ccc}
0&0 &0\\
0&0&e_j\\
0&0&0\end{array}\right),\quad R^{1,0}(e_k,\bar q)=
\left(\begin{array}{ccc}
0&-\bar e_k^t &0\\
0&0&e_k\\
0&0&0\end{array}\right),$

$R^{1,0}(q,\bar q)=\id_{\Co^{1,n+1}}-\frac{1}{2}\id_{\Co^m},$
where $1\leq j\leq m$, $m+1\leq k\leq n$.
\end{itemize}
\end{theorem}

In the case c), $\g\subset\u(1,1)$ is completely reducible, and it
corresponds to the symmetric space $\textrm{SL}(2,\Co)/\Co^*$ of
simple isometry group. Symmetric spaces corresponding to other
pairs have non-simple isometry groups and they are found in
\cite{K-O09}. Spaces corresponding to the cases a), b), d) and e)
are described also in \cite{Al10}.

\begin{corol} All indecomposable simply connected Calabi-Yau Lorentz-K\"ahler
symmetric spaces are exhausted by the cases a), b), d) and e) form
the above theorem. \end{corol}

This result is affirmative with the result from \cite{Al10}, where
it is shown that there are exactly two (up to isometry)
indecomposable simply connected symmetric Calabi-Yau pseudo-K\"ahler
manifolds in dimension 4, and the same holds for the dimension 6.

\bibliographystyle{unsrt}

\end{document}